\pdfoutput=1
\documentclass[10pt]{article}

\title{Order invariance of the generalized Rosenblatt transform on standard Borel spaces}
\author{Christian T.\ Covington \and Jeffrey W.\ Miller}
\date{\today}

\usepackage{setspace}
\usepackage[T1]{fontenc}
\usepackage[utf8]{inputenc}
\usepackage{lmodern}
\usepackage[margin=1in]{geometry}
\usepackage{microtype}

\usepackage{amsmath,amssymb,amsthm,mathtools}
\usepackage{aliascnt}
\usepackage{bm}
\usepackage{enumitem}
\newlist{questions}{enumerate}{1}
\setlist[questions]{label=\textup{(Q\arabic*)}, ref=\textup{Q\arabic*}}
\usepackage{booktabs,longtable,array}
\usepackage{graphicx}
\usepackage{capt-of}
\usepackage{xcolor}
\usepackage{tikz}
\usetikzlibrary{arrows.meta}

\usepackage[authoryear,round]{natbib}
\usepackage{etoolbox}
\AtBeginEnvironment{thebibliography}{\singlespacing}
\usepackage[hidelinks]{hyperref}
\hypersetup{
  pdftitle={\csname @title\endcsname},
  pdfauthor={Christian T. Covington and Jeffrey W. Miller}
}
\usepackage[nameinlink,capitalise,noabbrev]{cleveref}

\allowdisplaybreaks
\setlist[enumerate]{topsep=3pt,itemsep=2pt,parsep=0pt}

\theoremstyle{plain}
\newtheorem{theorem}{Theorem}[section]
\newaliascnt{proposition}{theorem}
\newtheorem{proposition}[proposition]{Proposition}
\aliascntresetthe{proposition}
\newaliascnt{lemma}{theorem}
\newtheorem{lemma}[lemma]{Lemma}
\aliascntresetthe{lemma}
\newaliascnt{corollary}{theorem}
\newtheorem{corollary}[corollary]{Corollary}
\aliascntresetthe{corollary}

\theoremstyle{definition}
\newaliascnt{definition}{theorem}
\newtheorem{definition}[definition]{Definition}
\aliascntresetthe{definition}
\newaliascnt{example}{theorem}
\newtheorem{example}[example]{Example}
\aliascntresetthe{example}

\theoremstyle{remark}
\newaliascnt{remark}{theorem}

\aliascntresetthe{remark}

\newcommand{\R}{\mathbb R}
\newcommand{\N}{\mathbb N}
\newcommand{\Q}{\mathbb Q}

\newcommand{\cA}{\mathcal A}
\newcommand{\cB}{\mathcal B}

\newcommand{\cE}{\mathcal E}

\newcommand{\cH}{\mathcal H}
\newcommand{\cI}{\mathcal I}

\newcommand{\cL}{\mathcal L}

\newcommand{\cS}{\mathcal S}
\newcommand{\cT}{\mathcal T}
\newcommand{\cX}{\mathcal X}

\newcommand{\fO}{\mathfrak{O}}

\newcommand{\set}[1]{\left\{#1\right\}}

\newcommand{\indep}{\mathrel{\perp\!\!\!\perp}}
\newcommand{\Law}{\mathcal L}
\newcommand{\Unif}{\operatorname{Unif}}

\crefname{theorem}{Theorem}{Theorems}
\Crefname{theorem}{Theorem}{Theorems}

\crefname{proposition}{Proposition}{Propositions}
\Crefname{proposition}{Proposition}{Propositions}

\crefname{lemma}{Lemma}{Lemmas}
\Crefname{lemma}{Lemma}{Lemmas}

\crefname{corollary}{Corollary}{Corollaries}
\Crefname{corollary}{Corollary}{Corollaries}

\crefname{definition}{Definition}{Definitions}
\Crefname{definition}{Definition}{Definitions}

\crefname{example}{Example}{Examples}
\Crefname{example}{Example}{Examples}

\crefname{remark}{Remark}{Remarks}
\Crefname{remark}{Remark}{Remarks}

\crefname{questionsi}{question}{questions}
\Crefname{questionsi}{Question}{Questions}

\begin{document}

\maketitle
\begin{abstract}
The Rosenblatt transform maps a random vector to independent uniform
variables through successive conditional distribution functions; for
laws on $\R^{\N}$ with atomless conditional distributions and
a fixed order of the coordinates, it is the almost surely unique
monotone triangular map with this property.
For a fixed reference law, transform-based tests remain valid under every fixed admissible order, but that choice can still affect their statistical power.
We extend the generalized Rosenblatt
transform to laws on countable products of standard Borel spaces, which
requires, in addition to an order on the index set, a total order on
each component space. For fixed orders, we establish identification and
uniqueness results for this extended generalized Rosenblatt transform (EGRT) that parallel those for
the classical transform. We then characterize order dependence in two senses: first demonstrating conditions under which the
EGRT is invariant under linear extensions of prescribed partial orders, and second showing that under every admissible order the output of the EGRT
retains all the statistical information from the input observations over reference-dominated families.
We close with a brief discussion of what these results mean for the relationship between order choice and goodness-of-fit testing.
\end{abstract}

\section{Introduction}
\label{sec:introduction}

The Rosenblatt transform~\citep{rosenblatt1952remarks} turns a random
vector with an absolutely continuous law on $\R^n$ into independent
standard uniform variables through iterative use of the
vector's conditional CDFs. When the conditional laws have atoms, the
conditional CDFs jump and the output is no longer uniform. Filling each
jump with an independent uniform variable allows the transform
to achieve uniformity for
arbitrary laws on $\R^n$
\citep{ruschendorf1981stochastically,brockwell2007universal}.
We call this randomized construction the \emph{generalized} Rosenblatt transform (GRT).

\citet{rosenblatt1952remarks} originally proposed the transform as a preprocessing
step for goodness-of-fit testing. If the conditional CDFs are taken
from a hypothesized model, which we call the \emph{reference} law, the
output is uniform on $[0,1]^n$ whenever the data follow the reference,
so any test of joint uniformity becomes a test of the model. This idea
underlies a range of multivariate diagnostics, historically developed for
real-valued coordinates (\cref{subsec:related_work}). Many statistical
models, however, describe data or parameters that are not real
vectors, e.g.\ graphs, partitions, matrices, or other structured objects.
Applying the transform to such a model requires a conditional CDF on
each component space, and hence a total order on that space, which we
call a \emph{component order}. On $\R$ the usual $\leq$ order serves
naturally and is generally ignored; on a space of graphs, for example,
there is no such default.

We formulate an extension of the generalized Rosenblatt transform on countable
products of standard Borel spaces, making explicit both the
component-space orders (our new addition) and the index order, which determines
which coordinates enter each conditioning set and has been noted as a relevant
design choice for the transform~\citep{rosenblatt1952remarks,berg2009copula}.
We call this formulation the \emph{extended generalized Rosenblatt transform}
(EGRT) and study, in three parts, how its properties depend on these order choices.

\begin{enumerate}[label=(\roman*)]
    \item \emph{Defining and constructing the EGRT (\cref{sec:egrt}):}
    We construct the EGRT for fixed admissible orders and show that
    its output is product-uniform precisely when the input law
    equals the reference. Any departure from the reference changes
    the law of some finite output prefix. We also establish
    almost-sure uniqueness within a class of locally randomized
    monotone triangular maps, extending the uniqueness result for the classical Rosenblatt transform~\citep{bogachev2005triangular}.
    \item \emph{Order invariance (\cref{sec:order_invariance}):}
    Given partial orders that encode the comparisons regarded as fixed,
    such as inclusion on a graph-valued component or the reachability
    order of a DAG on the index set, we ask when all of their
    \emph{admissible} extensions (the total orders that extend them and
    satisfy our requirements for a valid EGRT) yield the same EGRT.

    For a fixed index order, the key quantity is the conditional
    reference mass assigned to values incomparable with the
    observation. We call this \emph{conditional ambiguity}. If it
    vanishes almost everywhere, the transform is invariant to the
    choice of component order. The converse holds under a countable
    CDF-separation condition on the admissible extensions. We also
    show that CDF separation for every law is equivalent to the
    existence of a countable Borel representation of the component
    partial order.

    With the component orders fixed, index-order invariance is
    equivalent to \emph{ancestral screening}: each coordinate is
    conditionally independent of the coordinates at incomparable
    indices given those at its ancestors. Here we assume that the
    index partial order has finite ancestor sets. In particular,
    if the reference is Markov with respect to a DAG with finite
    ancestor sets, all of its finite-predecessor topological orders
    give the same transform. For any two EGRTs of the same reference,
    whatever their orders, agreement of the maps almost surely is
    equivalent to agreement of the output laws for every input law
    dominated by the reference.
    \item \emph{Statistical equivalence across orders (\cref{sec:statistical_equivalence}):}
    For each admissible order, a single measurable map recovers the
    observation from the output almost surely under every
    law whose finite-dimensional marginals are dominated by those
    of the reference. The original and transformed experiments
    are therefore equivalent over any such family, including every
    reference-dominated family.
    The behavior of a fixed test can vary between orders, but we show how
    one test can be transferred to its equivalent test under a different order.

    For reference-dominated laws, we characterize the output likelihood ratio as the input
    likelihood ratio evaluated at the reconstructed observation.
    This implies that every $f$-divergence between the output law and
    uniformity equals the corresponding divergence between the input
    law and the reference. For dominated alternatives, we also
    strengthen the finite-witness result in (i): each output prefix
    preserves the divergence between the corresponding input and
    reference marginals, and these divergences increase to the full
    divergence as the prefixes exhaust the index set.
\end{enumerate}

Proofs and supporting
results are presented in the supplementary material; references to
sections and results there use the prefix S.

\subsection{Related work}
\label{subsec:related_work}

Soon after~\citet{rosenblatt1952remarks}, \citet{knothe1957} introduced a closely related triangular map in the
study of convex bodies, and the resulting Knothe--Rosenblatt
rearrangement is now widely used as a transport map in sampling and
density estimation
\citep{baptista2024conditional,baptista2024representation}.

The Rosenblatt transform and its generalization were both
introduced on $\R^n$, but can be extended directly to countably many real
coordinates~\citep{bogachev2005triangular,ruschendorf1993regression}.
For components beyond the real line, related measurable representations
are known for Borel spaces: \citet{skorokhod1976representation} gives a
single conditional step, \citet{ruschendorf1993regression} note that
their sequential representation extends to processes with values in
Borel spaces, and \citet[Lemma~4.22]{kallenberg2021foundations} and
\citet[Theorem~A.3]{beiglboeck2025denseness} give jointly measurable
representations of probability kernels, the latter for atomless
conditional laws. In existing treatments of the Rosenblatt transform,
an order on a component space
appears, if at all, only implicitly, for example through a Borel
encoding into $\R$ \citep[Theorem~15.6]{kechris1995classical}.

Regarding uniqueness, \citet[Lemma~2.1]{bogachev2005triangular} show that, for laws on $\R^{\N}$
with atomless conditionals, the Rosenblatt transform is the unique uniformizing increasing triangular transform.
For laws with atoms, \citet[Proposition~5.8]{backhoff2017causal}
characterize the Knothe--Rosenblatt coupling on finite-dimensional
real spaces through bicausality and conditional monotonicity.
This determines the joint law of the observation and the output,
but does not specify how a given auxiliary uniform draw is used
to produce the output. Our uniqueness result concerns the
randomized map $(x,v)\mapsto T(x,v)$ itself.

In statistics, the transform is often used as a preprocessing step for goodness-of-fit
testing. Examples include multivariate
Kolmogorov--Smirnov tests \citep{justel1997multivariate},
goodness-of-fit tests for copulas
\citep{genest2009goodness,berg2009copula}, and the evaluation of
sequential and multivariate density forecasts
\citep{diebold1998,diebold1999multivariate,dovern2020order}.
Index-order effects on fixed diagnostics have long been recognized
\citep{rosenblatt1952remarks,ruschendorf2009distributional,berg2009copula}.
\citet{dovern2020order} construct order-invariant test statistics and
\citet{pournaderi2026calibrated} pools evidence across orders in attempts to
reduce the sensitivity of Rosenblatt-based diagnostics to the order.
Closest to our index-order results, \citet[Theorem~1]{akbari2023learning} show
that Knothe--Rosenblatt maps for compatible orders of a DAG agree up to
permutation in a finite, absolutely continuous structural equation
model with mechanisms strictly increasing in their independent noise
arguments and strictly increasing noise CDFs. Ancestral screening is
necessary as well as sufficient for index-order invariance, and its DAG
consequence extends this invariance to arbitrary standard Borel
components, local randomization, and countable index sets with finite
ancestor sets.

When introducing the idea of orders on the component space, we
start by assuming the existence of a partial order on that space.
Two order-theoretic questions will arise:
whether the partial order has an admissible extension, and whether its
admissible extensions can detect incomparable values. Both connect to
real-valued representations of preorders. An admissible extension
exists exactly when some countable family of monotone Borel functions
separates points, a Borel version of the injective monotones of
\citet[Lemmas~1--3]{hack2022injective}. Detecting incomparability
requires a countable family of Borel functions that represents the
order itself, so that $x$ precedes $y$ exactly when $f(x) \le f(y)$
for every $f$ in the family, called a countable multi-utility
\citep{hack2022injective,hack2023classification}.
\Cref{thm:countable_representation_cdf_separating} shows that, for a
partial order with Borel graph and a nonempty extension class, such a
representation exists exactly when the admissible extensions are
CDF-separating for every law.

\section{The extended generalized Rosenblatt transform}
\label{sec:egrt}

We begin by introducing standard results in our notation; further background can be found in
\citet{kechris1995classical} and
\citet[Chapter~1]{kallenberg2021foundations}.

A measurable space $(S,\cS)$ is a \emph{standard Borel space} if it is isomorphic to a Borel subset of a Polish space.
For measurable spaces $(S,\cS)$ and $(T,\cT)$, a \emph{probability kernel} from $S$ to $T$ is a map
\[
    K: S \times \cT \to [0,1]
\]
such that $K(s,\cdot)$ is a probability measure on $(T,\cT)$ for each $s \in S$, and $K(\cdot,B)$ is measurable for each $B \in \cT$.
For a measurable map $f:S \to T$ and a probability measure $\mu$ on $S$, we write $f_{\#}\mu$ for the pushforward measure,
defined by $(f_{\#}\mu)(B)=\mu(f^{-1}(B))$.

Let $\cI$ be a countable index set.
Throughout this section we fix a total order $\trianglelefteq$ on
$\cI$ in which every index has finitely many strict predecessors, so
that $(\cI,\trianglelefteq)$ is order-isomorphic to an initial segment
of $(\N,\leq)$.
We call such an order a
\emph{finite-predecessor index order}. We write $j \triangleleft k$ for
$j \trianglelefteq k$ with $j \neq k$, and call $D \subseteq \cI$
\emph{downward closed} if $d \in D$ and $i \trianglelefteq d$ imply
$i \in D$. A finite downward-closed set will also be called a
\emph{finite prefix}.

Let $(\cX_k,\cB_k)_{k \in \cI}$ be standard Borel spaces, with product space
\[
    \cX_{\cI}:=\prod_{k \in \cI}\cX_k,
    \qquad
    \cB_{\cI}:=\bigotimes_{k \in \cI}\cB_k.
\]
As a countable product of standard Borel spaces, $(\cX_{\cI},\cB_{\cI})$ is itself standard Borel~\citep{kechris1995classical,kallenberg2021foundations}.
For $S \subseteq \cI$, let $\pi_S:\cX_{\cI} \to \cX_S:=\prod_{j \in S}\cX_j$ be the coordinate projection and,
for a probability measure $\nu$ on $(\cX_{\cI},\cB_{\cI})$, let $\nu_S:=(\pi_S)_{\#}\nu$ denote the corresponding marginal.

For $k \in \cI$, let
\[
    P_{\triangleleft}(k):=\set{j \in \cI:j \triangleleft k}.
\]
For $x_{\cI} \in \cX_{\cI}$, write
\[
    x_{\triangleleft k}:=x_{P_{\triangleleft}(k)}=(x_j)_{j \triangleleft k},
    \qquad
    x_{\trianglelefteq k}:=(x_j)_{j \trianglelefteq k}=(x_{\triangleleft k},x_k),
\]
with corresponding product spaces $\cX_{\triangleleft k}:=\cX_{P_{\triangleleft}(k)}$ and $\cX_{\trianglelefteq k}$,
product $\sigma$-algebras $\cB_{\triangleleft k}$ and $\cB_{\trianglelefteq k}$,
projections $\pi_{\triangleleft k}$ and $\pi_{\trianglelefteq k}$, and
marginals $\nu_{\triangleleft k}$ and $\nu_{\trianglelefteq k}$.

To define the conditional CDF at coordinate $k$, we condition the
reference law on the preceding observations via the standard
disintegration theorem~\citep[Theorem~8.5]{kallenberg2021foundations}.

\begin{lemma}[Disintegration along the index order]
    \label{lem:disintegration_standard_borel}
    Let $\nu$ be a probability measure on $(\cX_{\cI},\cB_{\cI})$ and fix $k \in \cI$.
    There exists a \(\nu_{\triangleleft k}\)-a.e.\ unique probability kernel
    \[
    Q_k^\nu:\cX_{\triangleleft k} \times \cB_k \to [0,1]
    \]
    such that
    \[
    \nu\{x_{\cI}: x_{\triangleleft k} \in A,\ x_k \in B\}
    =
    \int_A Q_k^\nu(x_{\triangleleft k},B)\,\nu_{\triangleleft k}(dx_{\triangleleft k})
    \]
    for all \(A \in \cB_{\triangleleft k}\) and \(B \in \cB_k\).
    Moreover, $Q_k^\nu$ depends on the order $\trianglelefteq$ only through the strict predecessor set $P_{\triangleleft}(k)$.
\end{lemma}

We will also need conditional laws given other sets of coordinates.
For $S \subseteq \cI \setminus \set{k}$, write
\[
    Q_{k \mid S}^\nu:\cX_S \times \cB_k \to [0,1],
\]
for a fixed version of the conditional law of $X_k$ given $X_S$.
Such kernels exist by the same disintegration result; see also
\citet{blackwell1975,kallenberg2021foundations}.

An EGRT is defined pointwise using chosen kernel versions. Once $\nu$
and $\trianglelefteq$ are fixed, we choose one version of each
$Q_k^\nu$. Likewise, when another conditioning set $S$ is needed, we fix a
version of $Q_{k\mid S}^\nu$.
Let $\lambda$ denote Lebesgue measure on $[0,1]$ and
$\lambda^{\cI}$ its product.
For fixed orders, changing versions leaves the EGRT
unchanged $\nu\otimes\lambda^{\cI}$-a.e., and also
$\mu\otimes\lambda^{\cI}$-a.e.\ whenever $\mu\ll\nu$.

\subsection{Borel-embeddable component orders}
\label{subsec:borel_embeddable_orders}

The EGRT progresses by iteratively turning draws from each conditional law
$Q_k^\nu(x_{\triangleleft k},\cdot)$ into a uniform variable.
When defining the GRT on $\R$, the canonical total order $\leq$ is implicitly used to define CDFs.
On a general standard Borel space, however, an explicit total order is needed before one can define a CDF.
For our construction, we will use total orders that are Borel-embeddable in the following sense.

\begin{definition}[Borel-embeddable total order]
    \label{defn:borel_embeddable_total_order}
    Let $(\cX,\cB)$ be a standard Borel space.
    A total order $\preceq$ on $\cX$ is called a \emph{Borel-embeddable total order}
    if there exists a Borel injection $h_\preceq: \cX \to \R$ such that
    \[
        x \preceq y
        \qquad \Longleftrightarrow \qquad
        h_\preceq(x) \leq h_\preceq(y).
    \]
    We call any such $h_\preceq$ a \emph{Borel order embedding} for $\preceq$.
\end{definition}

An order embedding makes the order graph Borel, and hence makes its
lower sets measurable:
\[
    \set{(x,y):x \preceq y}=(h_\preceq \times h_\preceq)^{-1}\set{(s,t) \in \R^2:s \leq t}.
\]
Moreover, $h_\preceq$ identifies $(\cX,\preceq)$ with an ordered Borel
subset of $\R$, with a Borel inverse on its image
\citep[Corollary~15.2]{kechris1995classical}.
We write $x \prec y$ when $x \preceq y$ and $x \neq y$.

It is immediate that every standard Borel space admits such a total order by the Borel
isomorphism theorem \citep[Theorem~15.6]{kechris1995classical}.
However, when a component space has a canonical partial order, the analyst may want to use it,
and the isomorphism theorem does not guarantee that there exists a Borel-embeddable total order
which respects that partial order. We return to this question in~\cref{subsec:component_extensions}.

\subsection{The randomized probability integral transform}

The classical probability integral transform maps a real-valued random variable $X$ with distribution $\nu$ to $F_\nu(X)$,
but if $\nu$ has atoms, this ordinary transform is not uniform.
The randomized probability integral transform (rPIT)~\citep{ferguson1967mathematical} resolves this by
spreading the mass at each atom uniformly over the corresponding jump of the CDF.
Once a Borel-embeddable order has been chosen, this construction carries
over to any standard Borel space.

Fix a Borel-embeddable total order $\preceq$ on a standard Borel space $(\cX,\cB)$.
For a probability measure $\nu$ on $\cX$, define its CDF and strict-left CDF with respect to $\preceq$ by
\[
    F_\nu^\preceq(x):=\nu\{y \in \cX:y \preceq x\},
    \qquad
    F_\nu^\preceq(x^-):=\nu\{y \in \cX:y \prec x\}.
\]
These quantities are well-defined because Borel-embeddability makes the lower sets Borel.

\begin{definition}[Randomized probability integral transform]
    \label{defn:randomized_PIT}
    Let $\nu$ be a probability measure on a standard Borel space $(\cX,\cB)$ equipped with a Borel-embeddable total order $\preceq$.
    The \emph{randomized probability integral transform} with respect to $\nu$ and $\preceq$ is the map
    \[
        \mathrm{rPIT}_\nu^\preceq:\cX \times [0,1] \to [0,1],
        \qquad
        \mathrm{rPIT}_\nu^\preceq(x,v)
        :=F_\nu^\preceq(x^-)+v\,\nu(\set{x}).
    \]
\end{definition}

\cref{thm:randomized_PIT_uniform} reproduces, at our extended level of generality,
the classical result
that the rPIT is uniformizing~\citep[Proposition~2.1]{ruschendorf2009distributional}
and adds the converse, that $\nu$ is the only input law whose rPIT is uniform.

\begin{theorem}[Uniformity and identification for the rPIT]
    \label{thm:randomized_PIT_uniform}
    Let $(\cX,\cB)$ be a standard Borel space with a Borel-embeddable
    total order $\preceq$, let $\mu$ and $\nu$ be probability measures
    on $\cX$, and let $\lambda$ be Lebesgue measure on $[0,1]$. Then
    \[
        (\mathrm{rPIT}_\nu^\preceq)_{\#}(\mu\otimes\lambda)=\lambda
        \qquad\Longleftrightarrow\qquad
        \mu=\nu .
    \]
\end{theorem}

We also need this formula to be measurable in the
conditioning value. If $(S,\cS)$ is a measurable space and
$K:S\times\cB\to[0,1]$ is a probability kernel, the map
$(s,x,v)\mapsto\mathrm{rPIT}^{\preceq}_{K(s,\cdot)}(x,v)$ is jointly
measurable (\cref{lem:kernel_valued_rpit_measurability}). This allows
the one-coordinate construction to be combined sequentially to
create a triangular transform.

\subsection{Definition and uniformity of the EGRT}

We now have the necessary components to construct the EGRT.
Let $\lambda^{\cI}$ denote product Lebesgue measure on $[0,1]^{\cI}$.
After choosing versions of the kernels $(Q_k^\nu)_{k\in\cI}$ from \cref{lem:disintegration_standard_borel}, and Borel-embeddable total orders $\vec{\preceq}=(\preceq_k)_{k\in\cI}$ on the component spaces, \cref{lem:kernel_valued_rpit_measurability} shows that the coordinate functions below are measurable.
Because $\cI$ is countable, they define a measurable map
\[
    T_{\cI,\nu}^{\vec{\preceq},\trianglelefteq}:\cX_{\cI}\times[0,1]^{\cI}\to[0,1]^{\cI}
\]
coordinatewise by
\[
    \bigl[T_{\cI,\nu}^{\vec{\preceq},\trianglelefteq}(x_{\cI},v_{\cI})\bigr]_k
    :=
    \mathrm{rPIT}_{Q_k^\nu(x_{\triangleleft k},\cdot)}^{\preceq_k}(x_k,v_k),
    \qquad k\in\cI.
\]
For a $\trianglelefteq$-minimal index $k$, the predecessor space is $\cX_{\varnothing}=\set{\varnothing}$, and the kernel is the marginal $\nu_{\set{k}}$.
This map $T_{\cI, \nu}^{\vec{\preceq}, \trianglelefteq}$ is exactly our extended generalized Rosenblatt transform.

\begin{definition}[Extended generalized Rosenblatt transform]
    \label{defn:generalized_rosenblatt_transform}
    Given a reference law $\nu$, an index set $\cI$, an index order $\trianglelefteq$, and component-space orders $\vec{\preceq}$,
    the \emph{extended generalized Rosenblatt transform} (EGRT) is the map $T_{\cI,\nu}^{\vec{\preceq}, \trianglelefteq}$ above.
    For random elements $X_{\cI}$ of $\cX_{\cI}$ and $V_{\cI}\sim\lambda^{\cI}$ independent of $X_{\cI}$, we write
    \[
        U^{\vec{\preceq}, \trianglelefteq}
        :=T_{\cI,\nu}^{\vec{\preceq}, \trianglelefteq}(X_{\cI},V_{\cI}),
    \]
    and $\Law_\mu(\cdot)$ for laws computed when $X_{\cI}\sim\mu$.
\end{definition}

When varying only one of the two order choices, we suppress the other
from the notation and write $U_k^{\vec{\preceq}}$ or
$U_k^{\trianglelefteq}$. We write simply $U_k$ when both orders are
understood. Like the Rosenblatt transform and the GRT, the EGRT is a uniformizing transform.

\begin{theorem}[Uniformity of the EGRT]
    \label{thm:grt_uniformity}
    For every probability measure $\nu$ on $(\cX_{\cI},\cB_{\cI})$ and
    Borel-embeddable total orders $\vec{\preceq}$ on the component spaces,
    \[
        (T_{\cI,\nu}^{\vec{\preceq}, \trianglelefteq})_{\#}(\nu \otimes \lambda^{\cI})
        =\lambda^{\cI}.
    \]
    Equivalently, if $X_{\cI}\sim\nu$, the coordinates of
    $U^{\vec{\preceq}, \trianglelefteq}$ are i.i.d.\ $\mathrm{Unif}(0,1)$.
\end{theorem}

The reference law and orders in these results are fixed. Uniformity
is not automatic if either is chosen using the observations being
tested. For example, with one component and $\nu=\Unif[0,1]$, the
usual and reversed orders give $X$ and $1-X$, respectively. Choosing
the order that gives the smaller output produces $U=\min(X,1-X)$,
for which $\Pr_\nu(U\leq\alpha)=2\alpha$ when $0<\alpha<1/2$.
Thus a test of size $\alpha$ under either fixed order need not retain
that size after the order is selected from the data.

\subsection{Identification and finite witnesses of misspecification}
\label{subsec:identification_finite_witnesses}

When using Rosenblatt transforms for goodness-of-fit testing~\citep{genest2009goodness,berg2009copula},
one specifies a reference model $\nu$, observes data generated from an unknown law $\mu$, and assesses whether the transformed observations behave like independent uniforms.
The uniformizing property in \cref{thm:grt_uniformity} provides the implication
\[
    \mu=\nu
    \quad\Longrightarrow\quad
    (T_{\cI,\nu}^{\vec{\preceq}, \trianglelefteq})_{\#}(\mu \otimes \lambda^{\cI})
    =\lambda^{\cI}.
\]
Observed nonuniformity in transformed samples therefore provides evidence against $\mu=\nu$, exactly what is required for goodness-of-fit testing.
It is potentially of interest to also have the converse: if $(T_{\cI,\nu}^{\vec{\preceq}, \trianglelefteq})_{\#}(\mu \otimes \lambda^{\cI}) = \lambda^{\cI}$, then $\mu=\nu$.
For the classical Rosenblatt transform in the atomless Euclidean setting,
almost-sure inversion under the reference gives this converse for
reference-dominated input laws. For arbitrary input laws, however,
inversion on a reference-full set alone does not suffice, since the
input may charge its complement. The result below imposes no
domination assumption; its proof uses one-coordinate identification
and proceeds by induction over finite prefixes.

\begin{theorem}[Identification by EGRT uniformity]
    \label{thm:grt_converse_uniformization}
    Let $\mu$ and $\nu$ be probability measures on $(\cX_{\cI},\cB_{\cI})$, and fix an index order $\trianglelefteq$ and component-space orders $\vec{\preceq}$.
    Then
    \[
        \bigl(T_{\cI,\nu}^{\vec{\preceq}, \trianglelefteq}\bigr)_{\#}(\mu \otimes \lambda^{\cI})=\lambda^{\cI}
        \qquad \Longleftrightarrow \qquad
        \mu=\nu.
    \]
\end{theorem}

In practice, goodness-of-fit tests will take in only finitely many elements at once, so we want to be sure
that evidence against the hypothesis $\mu = \nu$ manifests in finitely many coordinates.
\cref{cor:finite_witness_misspecification} guarantees this.

\begin{corollary}[Finite witness of misspecification]
    \label{cor:finite_witness_misspecification}
    In the setting of~\cref{thm:grt_converse_uniformization}, if $\mu \neq \nu$,
    then there exists a finite $\trianglelefteq$-downward closed set $D \subseteq \cI$ such that
    \[
        \bigl(\pi_D \circ T_{\cI,\nu}^{\vec{\preceq}, \trianglelefteq}\bigr)_{\#}
        (\mu \otimes \lambda^{\cI})
        \neq \lambda^D.
    \]
\end{corollary}

\subsection{Uniqueness of the EGRT}
\label{subsec:egrt_uniqueness}

A natural question is whether the EGRT is the only map satisfying~\cref{thm:grt_uniformity}.
For the classical Rosenblatt transform on $\R^\N$ and a reference measure $\nu$ with atomless conditional measures,
\citet{bogachev2005triangular} prove uniqueness after restricting to deterministic increasing triangular maps.
Our setting is different in that the EGRT acts on the enlarged space $\cX_{\cI}\times[0,1]^{\cI}$ and uses
auxiliary randomness.
The reconstruction associated with the randomized PIT recovers $X_{\cI}$ from the transformed coordinates almost surely,
but it does not generally recover the full auxiliary vector $V_{\cI}$; in particular, atomless conditional laws make the corresponding $V_k$ irrelevant.
Thus, our setting's uniqueness statement uses a restricted notion of triangularity: the $k$-th output may depend on $x_{\trianglelefteq k}$ and on the current randomizer $v_k$, but not on earlier randomizers $v_{\triangleleft k}$.

\begin{definition}[Locally randomized triangular map]
    \label{defn:triangular_map}
    A measurable map $T:\cX_{\cI} \times [0,1]^{\cI} \to [0,1]^{\cI}$ is \emph{locally randomized $\trianglelefteq$-triangular} if, for each $k \in \cI$, there exists a measurable function
    \[
        T_k:\cX_{\trianglelefteq k} \times [0,1] \to [0,1]
    \]
    such that
    \[
        [T(x_{\cI},v_{\cI})]_k=T_k(x_{\trianglelefteq k},v_k)
    \]
    for every $(x_{\cI},v_{\cI})$.
\end{definition}

We then define monotonicity with respect to both our observed data and the auxiliary uniform.

\begin{definition}[Monotonicity in the current pair]
    \label{defn:order_monotone_map}
    Let $\nu$ be a probability measure on $(\cX_{\cI},\cB_{\cI})$ and
    $\vec{\preceq}=(\preceq_k)_{k \in \cI}$ be Borel-embeddable total orders on the component spaces.
    A locally randomized triangular map $T$ is \emph{$\vec{\preceq}$-monotone relative to $\nu$} if, for every $k \in \cI$ and $\nu_{\triangleleft k}$-a.e.\ history $x_{\triangleleft k}$, the map
    \[
        (x_k,v_k) \longmapsto T_k(x_{\triangleleft k},x_k,v_k)
    \]
    is nondecreasing for the lexicographic order on $\cX_k \times [0,1]$; that is, whenever $x_k \prec_k x'_k$, or $x_k=x'_k$ and $v_k \leq v'_k$,
    \[
        T_k(x_{\triangleleft k},x_k,v_k)
        \leq
        T_k(x_{\triangleleft k},x'_k,v'_k).
    \]
\end{definition}

For fixed orders, these two structural conditions, together with uniformization, uniquely characterize the EGRT up to null sets.

\begin{theorem}[Uniqueness among locally randomized monotone triangular maps]
    \label{thm:grt_uniqueness}
    Fix a probability measure $\nu$ on $(\cX_{\cI},\cB_{\cI})$ and Borel-embeddable total orders $\vec{\preceq}$ on the component spaces.
    The extended generalized Rosenblatt transform $T_{\cI,\nu}^{\vec{\preceq}, \trianglelefteq}$ is the $\nu \otimes \lambda^{\cI}$-a.e.\ 
    unique measurable map $T:\cX_{\cI} \times [0,1]^{\cI} \to [0,1]^{\cI}$ satisfying:
    \begin{enumerate}[label=(\roman*)]
        \item $T$ is uniformizing, $T_{\#}(\nu \otimes \lambda^{\cI})=\lambda^{\cI}$;
        \item $T$ is locally randomized $\trianglelefteq$-triangular in the sense of \cref{defn:triangular_map};
        \item $T$ is $\vec{\preceq}$-monotone relative to $\nu$, in the sense of \cref{defn:order_monotone_map}.
    \end{enumerate}
\end{theorem}

\section{Invariance under choice of order}
\label{sec:order_invariance}

Every EGRT constructed from the same reference law $\nu$ has output law $\lambda^{\cI}$ under
$\nu\otimes\lambda^{\cI}$ by~\cref{thm:grt_uniformity}, a property we call \emph{reference invariance}.
For goodness-of-fit testing, this means that a fixed test of uniformity applied to the output has the same size under every
fixed admissible order.
We now introduce two other levels of invariance, which are of interest for goodness-of-fit applications.
At the strongest level, which we call \emph{a.e.\ map invariance}, two transforms agree
for $\nu\otimes\lambda^{\cI}$-almost all values of observed data and auxiliary uniform.
Thus, any fixed diagnostic calculated over these transforms will yield identical answers almost surely when the same
auxiliary uniforms are used.
Between these two levels sits \emph{distributional invariance}, under which the realized maps may differ in general,
but their output laws coincide under every input law dominated by the reference.
Every fixed test then has the same power against every reference-dominated alternative under both transforms,
although its realized values need not agree for general randomized maps.
Writing $T$ and $T'$ for two randomized maps,
\begin{equation}
    \label{eq:invariance_hierarchy}
    \begin{aligned}
        &\text{a.e.\ map invariance:}
        && T=T'\quad\nu \otimes \lambda^{\cI}\text{-a.e.}\\
        &\qquad\Downarrow\\
        &\text{distributional invariance:}
        && T_{\#}(\mu \otimes \lambda^{\cI})
           =T'_{\#}(\mu \otimes \lambda^{\cI})
           \quad\text{for all }\mu \ll \nu\\
        &\qquad\Downarrow\\
        &\text{reference invariance:}
        && T_{\#}(\nu \otimes \lambda^{\cI})
           =T'_{\#}(\nu \otimes \lambda^{\cI}).
    \end{aligned}
\end{equation}
The first implication holds because $\mu \ll \nu$ implies $\mu \otimes \lambda^{\cI} \ll \nu \otimes \lambda^{\cI}$,
and the second follows by taking $\mu=\nu$.
For EGRTs, the first implication reverses, so a.e.\ map and distributional invariance are in fact equivalent.
Each coordinate of an EGRT is a nondecreasing affine function of its own auxiliary uniform,
so its law given the observation determines it,
and conditioning the reference on a set where two transforms differ separates their output laws.

To compare maps against different index orders, we evaluate them at a
fixed auxiliary vector $V$ and index all outputs by their fixed
coordinate labels $k \in \cI$, rather than by their position
under the order. Thus, reordering $\cI$ alters conditioning sets
but preserves the identity of each output coordinate.
As such, a diagnostic based on the first $m$ positions of
the chosen order is not order-invariant, as permuting the indices
changes which subset of $\cI$ is used.

\begin{lemma}[Map and distributional invariance coincide]
    \label{lem:map_distributional_equivalence}
    Let $\nu$ be a probability measure on $(\cX_{\cI},\cB_{\cI})$, and let
    $T:=T_{\cI,\nu}^{\vec{\preceq},\trianglelefteq}$ and
    $T':=T_{\cI,\nu}^{\vec{\preceq}',\trianglelefteq'}$ be EGRTs of $\nu$
    built from any finite-predecessor index orders $\trianglelefteq,\trianglelefteq'$
    and any Borel-embeddable component orders $\vec{\preceq},\vec{\preceq}'$.
    Then $T=T'$ $\nu\otimes\lambda^{\cI}$-a.e.\ if and only if
    $T_{\#}(\mu\otimes\lambda^{\cI})=T'_{\#}(\mu\otimes\lambda^{\cI})$
    for every $\mu\ll\nu$.
    In that case, $T=T'$ also holds $\mu\otimes\lambda^{\cI}$-a.e.\ for every $\mu\ll\nu$.
\end{lemma}

The two main results of this section, \cref{thm:component_order_exact_invariance,thm:exact_index_order_invariance},
characterize when these two equivalent forms hold.

We ask when an EGRT is essentially invariant to the choice of orders used to construct it.
In goodness-of-fit testing, these orders must be chosen before any test can be run,
and on spaces such as graphs there is rarely a canonical total order;
invariance identifies when the choice affects neither the power of a test nor its outcome.
The two order choices affect the transform in different ways, so we study their invariance properties separately.
For invariance with respect to the component-space orders, we always start from a partial order on each component space,
representing the comparisons the analyst regards as fixed, and the relevant decision is how to choose a total order consistent with it.
Many natural spaces of interest come endowed with canonical partial orders,
e.g.\ edge inclusion on the space of graphs on a fixed vertex set,
 or the L\"owner order on the cone of symmetric positive definite matrices.
When no such order exists and the analyst is unwilling to choose one,
the trivial partial order, under which distinct elements are incomparable,
still falls under our theory.

\subsection{Borel-embeddable extensions of component-space partial orders}
\label{subsec:component_extensions}

A binary relation $\preceq_p$ is a \emph{partial order} on a set $\cX$ if it is reflexive, antisymmetric, and transitive.
A total order $\preceq$ is a \emph{linear extension} of $\preceq_p$ if $x \preceq_p y$ implies $x \preceq y$.
For a partial order on a standard Borel space, write
\[
    \cE(\preceq_p)
    :=
    \set{\preceq:
    \preceq \text{ is a Borel-embeddable linear extension of } \preceq_p}.
\]
We call the elements of $\cE(\preceq_p)$ \emph{admissible} extensions of $\preceq_p$.

We require two properties of this class. The first is nonemptiness,
so that the component order can be chosen; see
\cref{prop:countable_extension,thm:admissibility_characterization}.
The second is that whenever two independent draws
from a law are incomparable with positive probability, some two
admissible extensions give different strict-left CDFs on a set of
positive measure under that law. This second property is needed for
the converse from invariance to zero conditional ambiguity;
nonemptiness alone only ensures that the transform can be constructed.
A countable version of this property
is characterized in \cref{thm:countable_representation_cdf_separating}
and used in \cref{thm:component_order_exact_invariance}.

Szpilrajn's extension theorem guarantees an abstract linear extension~\citep{szpilrajn1930extension},
but it does not on its own provide the Borel order embedding required by the EGRT.
For countable spaces, Borel-embeddability of linear extensions is automatic.

\begin{proposition}[Borel-embeddability on countable spaces]
    \label{prop:countable_extension}
    Let $(\cX,\cB)$ be a countable standard Borel space and let $\preceq_p$ be a partial order on $\cX$.
    Every linear extension of $\preceq_p$ is a Borel-embeddable total order on $\cX$.
\end{proposition}

For general standard Borel spaces, admissible extensions correspond to Borel injections into $\R$ that are strictly increasing along the partial order.
Given such an injection $h$, the pullback relation $x\preceq_h y:\Leftrightarrow h(x)\leq h(y)$ is an admissible extension of $\preceq_p$,
and conversely a Borel order embedding for any member of $\cE(\preceq_p)$ is such an injection.
Such an injection can in turn be built by combining any countable family of Borel functions
that are monotone along $\preceq_p$ and together separate points.
\citet[Lemmas~1--3]{hack2022injective} construct injective monotones in this way from countable separating families of increasing sets;
see also \citet[Lemma~1(iii)]{hack2023classification} for the corresponding characterization.
Together, these observations give a criterion for the extension class to be nonempty.

\begin{theorem}[Borel-embeddable extension criterion]
    \label{thm:admissibility_characterization}
    Let $(\cX,\cB)$ be a standard Borel space and let $\preceq_p$ be a
    partial order on $\cX$. The following are equivalent:
    \begin{enumerate}[label=(\roman*)]
        \item $\cE(\preceq_p)\neq\varnothing$;
        \item there is a Borel injection $h:\cX\to\R$ such that
              $x\prec_p y$ implies $h(x)<h(y)$;
        \item there is a countable family of Borel functions
              $f_n:\cX\to\R$ that separates points and is monotone
              with respect to $\preceq_p$, in the sense that
              $x\preceq_p y$ implies $f_n(x)\leq f_n(y)$ for every $n$.
    \end{enumerate}
    Given a family in (iii), enumerate $\N\times\Q$ as
    $((n_j,r_j))_{j\geq1}$ and set
    \[
        h(x):=\sum_{j=1}^\infty3^{-j}\mathbf1\{f_{n_j}(x)>r_j\}.
    \]
    This $h$ satisfies (ii), and its pullback order belongs to
    $\cE(\preceq_p)$.
\end{theorem}

We now move to the second question.
The only values that can change sides of an observation between total orders
are those left incomparable by the partial order.

\begin{definition}[Incomparability set of a point]
    \label{defn:incomparable_set}
    Let $(\cX,\cB)$ be a standard Borel space and let $\preceq_p$ be a partial
    order on $\cX$.  For $x \in \cX$, define
    \[
        I_p(x)
        :=
        \set{y \in \cX:
        y \npreceq_px\text{ and }x \npreceq_py}.
    \]
\end{definition}

\begin{definition}[CDF-separating extension class]
    \label{defn:cdf_separating_extension_class}
    Let $\preceq_p$ be a partial order with Borel graph on a standard Borel
    space $(\cX,\cB)$ and let $\cH$ be a family of probability measures on
    $\cX$.  We call $\cE(\preceq_p)$ \emph{CDF-separating for $\cH$} if it
    contains a countable subfamily $\cE_0(\preceq_p)$ such that, for every
    $\eta \in \cH$ with $\eta\set{x:\eta(I_p(x))>0}>0$, there exist
    $\preceq, \preceq' \in \cE_0(\preceq_p)$ with
    \[
        \eta\set{x:
        F_\eta^\preceq(x^-)
        \neq
        F_\eta^{\preceq '}(x^-)}>0.
    \]
    We say simply that $\cE(\preceq_p)$ is \emph{CDF-separating} if this holds
    for the family of all probability measures on $\cX$.
\end{definition}

We require the same countable subfamily for all laws in $\cH$.

The answer to the second question strengthens the condition of \cref{thm:admissibility_characterization}.
Instead of only separating points, the statistics must determine the partial order, in the sense that
$x\preceq_p y$ exactly when $f_n(x)\leq f_n(y)$ for every $n$.
Such families are the countable multi-utilities of \citet{hack2022injective,hack2023classification},
and \citet[Proposition~6]{hack2022injective} show that a countable multi-utility can be replaced by one consisting of injective monotones.

The functions $h_{n,q}$ in \cref{thm:countable_representation_cdf_separating} are a Borel version of this replacement.
Each is a Borel injective monotone, so its pullback order is admissible by \cref{thm:admissibility_characterization},
and for every incomparable pair, the family contains orders ranking them in both directions.
\Cref{thm:countable_representation_cdf_separating} shows that this family is CDF-separating for every law, including laws without atoms, and that, for a partial order with a nonempty extension class, CDF separation for every law in turn requires a countable multi-utility.
Countable spaces are a special case; by reflexivity and transitivity, the indicators $\mathbf 1\{a\preceq_p\cdot\}$, $a\in\cX$, form a countable multi-utility (see also \cref{lem:countable_cdf_separating}).

\begin{theorem}[CDF separation and countable order representations]
    \label{thm:countable_representation_cdf_separating}
    Let $(\cX,\cB)$ be a standard Borel space with partial order
    $\preceq_p$. The following are equivalent:
    \begin{enumerate}[label=(\roman*)]
        \item there are Borel functions $f_n:\cX\to\R$, $n\geq1$,
        such that
        \begin{equation}
            \label{eq:countable_order_representation}
            x\preceq_p y
            \quad\Longleftrightarrow\quad
            f_n(x)\leq f_n(y)\ \text{for every }n\geq1;
        \end{equation}
        \item $\preceq_p$ has Borel graph and $\cE(\preceq_p)$ is
        nonempty and CDF-separating for all probability measures on
        $\cX$.
    \end{enumerate}

    Given a family in (i), choose a Borel injection $h_0:\cX\to(0,1)$
    that is strictly increasing along $\prec_p$, as supplied by
    \cref{thm:admissibility_characterization}. A countable separating
    family consists of the pullback orders of
    \[
        h_{n,q}(x):=\mathbf1\{f_n(x)>q\}+h_0(x),
        \qquad n\geq1,\ q\in\Q.
    \]
\end{theorem}

\begin{example}[Graph inclusion]
    \label{ex:graph_inclusion_extensions}
    Let $\cX=\{0,1\}^{E}$ be the space of simple graphs on a fixed
    countably infinite labeled vertex set, where $E$ is the countable
    set of possible edges. Equip $\cX$ with its product Borel
    $\sigma$-algebra and order it by edge inclusion. The coordinate
    functions $f_e(G):=\mathbf1\{e\in G\}$ are Borel and satisfy
    \[
        G\subseteq H
        \quad\Longleftrightarrow\quad
        f_e(G)\leq f_e(H)\ \text{for every }e\in E.
    \]
    By \cref{thm:countable_representation_cdf_separating}, the class
    of admissible extensions is therefore nonempty and CDF-separating
    for every probability law on this uncountable graph space.
\end{example}

\subsection{Invariance across component-space extensions}
\label{subsec:component_exact_invariance}

\Cref{subsec:component_extensions} establishes when the class of admissible extensions of a partial order is nonempty, and when it is CDF-separating.
An admissible extension is typically not unique, so it is of interest to know how much the choice matters.
We start by setting up the building block of this analysis, \emph{conditional ambiguity}.

Fix a finite-predecessor index order $\trianglelefteq$ and the reference kernels $(Q_k^\nu)_{k \in \cI}$.
For each $k$, suppose that $\cX_k$ carries a partial order $\preceq_{p,k}$ with Borel graph and a
nonempty class $\cE(\preceq_{p,k})$ of admissible extensions.
Write $I_{p,k}(x)$ for the incomparability set of~\cref{defn:incomparable_set}
computed in $(\cX_k, \preceq_{p,k})$.

Two admissible extensions of $\preceq_{p,k}$ can disagree about whether a point lies
below or above $x_k$ only if that point is incomparable to $x_k$ and thus in $I_{p,k}(x_k)$.
Therefore, given the predecessor history $x_{\triangleleft k}$, the reference conditional probability of this set,
which we call the \emph{conditional ambiguity},
bounds how much the choice of extension can matter at $x_k$.

\begin{definition}[Conditional ambiguity]
    \label{defn:conditional_ambiguity}
    For a probability measure $\nu$, the \emph{conditional ambiguity} at coordinate $k$ is
    \[
        \Delta_k(x_{\trianglelefteq k})
        :=
        Q_k^\nu\bigl(x_{\triangleleft k},I_{p,k}(x_k)\bigr).
    \]
    This map is $\cB_{\trianglelefteq k}$-measurable, so a.e.\ statements about it are well posed; see \cref{lem:conditional_ambiguity_measurable}.
\end{definition}

For a graph-valued component ordered by inclusion as in
\cref{ex:graph_inclusion_extensions}, $\Delta_k$ is the reference
conditional probability of a graph that neither contains the observed
graph nor is contained in it.

Because conditional ambiguity bounds how far the choice of extension can move each coordinate of the transform,
ensuring that the ambiguity is 0 forces every admissible choice to produce the same transform.

\begin{proposition}[Ambiguity bounds the effect of the component order]
    \label{prop:zero_ambiguity_sufficient}
    For every pair of extension sequences
    $\vec{\preceq},\vec{\preceq}' \in \prod_{k \in \cI}\cE(\preceq_{p,k})$,
    every $k\in\cI$, and every $(x,v)\in\cX_{\cI}\times[0,1]^{\cI}$,
    \[
        \bigl|[T_{\cI,\nu}^{\vec{\preceq}}(x,v)]_k
        -[T_{\cI,\nu}^{\vec{\preceq}'}(x,v)]_k\bigr|
        \leq \Delta_k(x_{\trianglelefteq k}).
    \]
    Consequently, if $\Delta_k=0$ $\nu_{\trianglelefteq k}$-a.e.\ for
    every $k$, then $T_{\cI,\nu}^{\vec{\preceq}}=T_{\cI,\nu}^{\vec{\preceq}'}$
    $\mu\otimes\lambda^{\cI}$-a.e.\ for every $\mu\ll\nu$, including
    $\mu=\nu$, and the two transforms have the same output law under
    every such $\mu$.
\end{proposition}

\begin{example}[Ambiguity on a space of graphs]
    \label{ex:graph_component_ambiguity}
    Consider a single graph-valued component on a fixed vertex set,
    with two possible edges $e_1,e_2$. Identifying each graph with its
    edge set, the component space consists of
    \[
        H_0=\varnothing,\qquad H_1=\{e_1\},\qquad
        H_2=\{e_2\},\qquad H_{12}=\{e_1,e_2\}.
    \]
    Order these graphs by inclusion and let $\nu$ be uniform on the
    four graphs. The only incomparable pair is $(H_1,H_2)$, so the
    ambiguity is $\Delta(H_1)=\Delta(H_2)=1/4$, while
    $\Delta(H_0)=\Delta(H_{12})=0$. The two admissible extensions are
    \[
        \begin{aligned}
            \preceq &: H_0\prec H_1\prec H_2\prec H_{12},\\
            \preceq' &: H_0\prec' H_2\prec' H_1\prec' H_{12}.
        \end{aligned}
    \]
    At the observation $H_1$, their transforms are
    \[
        \mathrm{rPIT}^{\preceq}_{\nu}(H_1,v)=\frac{1+v}{4},
        \qquad
        \mathrm{rPIT}^{\preceq'}_{\nu}(H_1,v)=\frac{2+v}{4}.
    \]
    Changing the order moves the transformed value by exactly $1/4$,
    attaining the ambiguity bound in
    \cref{prop:zero_ambiguity_sufficient}. The transform at $H_2$
    moves by the same amount in the opposite direction, while the
    transforms at $H_0$ and $H_{12}$ are unchanged.

    The upper-tail test $\psi(u)=\mathbf1\{u>1/2\}$ has size
    $1/2$ under both orders. Under $\preceq$ it rejects exactly when
    $e_2$ is present, whereas under $\preceq'$ it rejects exactly
    when $e_1$ is present, almost surely. Thus, under
    $\mu=\nu(\cdot\mid e_1\text{ is present})$, its powers are
    $1/2$ and $1$, respectively. We return to this example in
    \cref{ex:order_test_transfer} to show how the test can be
    transferred between orders.
\end{example}

For the converse, that invariance implies zero ambiguity, the class of extensions must be rich enough to reveal positive ambiguity,
which is why we need~\cref{defn:cdf_separating_extension_class}.
At coordinate $k$ the relevant family is
\begin{equation}
    \label{eq:reference_kernel_family}
    \cH_k^\nu
    :=
    \set{Q_k^\nu(x_{\triangleleft k},\cdot):
    x_{\triangleleft k} \in \cX_{\triangleleft k}},
\end{equation}
the conditional laws that the EGRT actually uses at $k$.
The class $\cE(\preceq_{p,k})$ is CDF-separating for $\cH_k^\nu$, whatever the reference kernels, when $\cX_k$ is countable (\cref{lem:countable_cdf_separating}) or,
more generally, satisfies~\cref{eq:countable_order_representation} (\cref{thm:countable_representation_cdf_separating}).

\begin{theorem}[Component-order invariance]
    \label{thm:component_order_exact_invariance}
    Assume that each $\preceq_{p,k}$ has Borel graph and
    $\cE(\preceq_{p,k})\neq\varnothing$, and consider the
    following conditions:
    \begin{enumerate}[label=(\roman*)]
        \item \textbf{zero conditional ambiguity a.e.:} for every $k \in \cI$,
        \[
            \Delta_k=0
            \qquad \nu_{\trianglelefteq k}\text{-a.e.};
        \]
        \item \textbf{a.e.\ map invariance:} for every pair of extension sequences
        $\vec{\preceq},\vec{\preceq}' \in
        \prod_{k \in \cI}\cE(\preceq_{p,k})$,
        \[
            T_{\cI,\nu}^{\vec{\preceq}}(x,v)
            =
            T_{\cI,\nu}^{\vec{\preceq}'}(x,v)
            \qquad
            \nu \otimes \lambda^{\cI}\text{-a.e.};
        \]
        \item \textbf{distributional invariance under reference-dominated measures:} for every $\mu \ll \nu$ and pair of extension sequences,
        \[
            (T_{\cI,\nu}^{\vec{\preceq}})_{\#}
            (\mu \otimes \lambda^{\cI})
            =
            (T_{\cI,\nu}^{\vec{\preceq}'})_{\#}
            (\mu \otimes \lambda^{\cI}).
        \]
    \end{enumerate}
    Then \textup{(i)}$\Rightarrow$\textup{(ii)}$\Leftrightarrow$\textup{(iii)}.
    If, in addition, $\cE(\preceq_{p,k})$ is CDF-separating for $\cH_k^\nu$
    of~\cref{eq:reference_kernel_family} for every $k \in \cI$, then
    \textup{(ii)}$\Rightarrow$\textup{(i)}, so all three conditions are equivalent.
    Under any of these conditions, the equality in~(ii) also holds $\mu \otimes \lambda^{\cI}$-a.e.\ for every $\mu \ll \nu$.
\end{theorem}

\Cref{prop:zero_ambiguity_sufficient} gives (i)$\Rightarrow$(ii), and \cref{lem:map_distributional_equivalence} gives (ii)$\Leftrightarrow$(iii); neither requires CDF separation, which is used only for (ii)$\Rightarrow$(i).

\subsection{Invariance across index-order extensions}
\label{subsec:index_exact_invariance}

We now fix Borel-embeddable component-space orders $\vec{\preceq}=(\preceq_k)_{k\in\cI}$ and vary the order on the index set.
A component-order extension changes which values of a fixed conditional distribution count as lower than the observation, while an index order changes the conditioning set itself: under a finite-predecessor index order $\trianglelefteq$, the kernel at coordinate $k$ is $Q_k^\nu=Q_{k\mid P_{\triangleleft}(k)}^\nu$.
Exact invariance must consequently come from a probabilistic reduction of these conditional kernels.
We write $T_{\cI,\nu}^{\trianglelefteq}$ when the fixed component orders are suppressed.

Let $\trianglelefteq_p$ be a partial order on $\cI$.  For $k \in \cI$,
define its strict ancestors and incomparable indices by
\[
    \mathrm{anc}_p(k)
    :=
    \set{j:j \triangleleft_p k},
    \qquad
    \mathrm{inc}_p(k)
    :=
    \set{j \in \cI \setminus \set{k}:
    j \ntrianglelefteq_p k\text{ and }k \ntrianglelefteq_p j}.
\]
Let $\cL(\trianglelefteq_p)$ be the set of finite-predecessor linear
extensions of $\trianglelefteq_p$; we call its elements \emph{admissible}
index orders.  If every ancestor set is finite,
this class is nonempty (\cref{lem:index_linear_extensions_exist}).
For every $\trianglelefteq \in \cL(\trianglelefteq_p)$,
\[
    \mathrm{anc}_p(k)
    \subseteq
    P_{\triangleleft}(k)
    \subseteq
    \mathrm{anc}_p(k) \cup \mathrm{inc}_p(k).
\]

We now have what we need for the index-order counterpart of~\cref{thm:component_order_exact_invariance}.
An admissible index order cannot change whether an ancestor of $k$ precedes $k$: every element of $\mathrm{anc}_p(k)$ belongs to $P_{\triangleleft}(k)$ under every $\trianglelefteq\in\cL(\trianglelefteq_p)$.
What may change is which finite subset of the indices incomparable with $k$ is also placed before $k$.

\begin{theorem}[Index-order invariance]
    \label{thm:exact_index_order_invariance}
    Let $\trianglelefteq_p$ be a partial order on the countable set $\cI$, with finite $\mathrm{anc}_p(k)$ for every $k \in \cI$,
    and $\nu$ be a probability measure on $(\cX_{\cI},\cB_{\cI})$.
    The following are equivalent:
    \begin{enumerate}[label=(\roman*)]
        \item \textbf{ancestral screening of incomparable indices:} for every $k \in \cI$,
        \[
            X_k
            \indep_\nu
            X_{\mathrm{inc}_p(k)}
            \mid
            X_{\mathrm{anc}_p(k)};
        \]
        \item \textbf{ancestor-kernel reduction:} for every $k \in \cI$ and every $\trianglelefteq \in \mathcal L(\trianglelefteq_p)$,
        \[
            Q_{k \mid P_{\triangleleft}(k)}^\nu
            \bigl(X_{P_{\triangleleft}(k)},\cdot\bigr)
            =
            Q_{k \mid \mathrm{anc}_p(k)}^\nu
            \bigl(X_{\mathrm{anc}_p(k)},\cdot\bigr)
            \qquad \nu\text{-a.e.};
        \]
        \item \textbf{a.e.\ map invariance:} for every pair
        $\trianglelefteq, \trianglelefteq ' \in \mathcal L(\trianglelefteq_p)$,
        \[
            T_{\cI,\nu}^{\trianglelefteq}(x,v)
            =
            T_{\cI,\nu}^{\trianglelefteq '}(x,v)
            \qquad
            \nu \otimes \lambda^{\cI}\text{-a.e.};
        \]
        \item \textbf{distributional invariance under reference-dominated measures:} for every $\mu \ll \nu$ and every pair
        $\trianglelefteq, \trianglelefteq ' \in \mathcal L(\trianglelefteq_p)$,
        \[
            \bigl(T_{\cI,\nu}^{\trianglelefteq}\bigr)_{\#}
            (\mu \otimes \lambda^{\cI})
            =
            \bigl(T_{\cI,\nu}^{\trianglelefteq '}\bigr)_{\#}
            (\mu \otimes \lambda^{\cI}).
        \]
    \end{enumerate}
    Under any of these conditions, the equality in~\textup{(iii)} also holds
    $\mu \otimes \lambda^{\cI}$-a.e.\ for every $\mu \ll \nu$.
\end{theorem}

As in \cref{thm:component_order_exact_invariance}, the equivalence (iii)$\Leftrightarrow$(iv) is an instance of \cref{lem:map_distributional_equivalence};
the content specific to index orders is (i)$\Leftrightarrow$(ii)$\Leftrightarrow$(iii).

Ancestral screening has a familiar interpretation through the local Markov
property for Bayesian networks~\citep{lauritzen1996graphical}.
Let $G=(\cI,E)$ be a DAG on a countable vertex set, and let its \emph{reachability relation} be
\[
    j \trianglelefteq_pk
    \quad: \Longleftrightarrow \quad
    j=k\text{, or there is a directed path from }j\text{ to }k\text{ in }G.
\]
Acyclicity makes this a partial order such that its ancestors are the ordinary graph ancestors, and its incomparable indices are the non-ancestor non-descendants.
Write $\mathrm{Pa}_G(k):=\set{j:(j,k) \in E}$ for the parents of $k$, call $j$ a
\emph{descendant} of $k$ when $k \triangleleft_pj$, and write
$\mathrm{ND}_G(k):=\cI \setminus (\set k \cup \set{j:k \triangleleft_pj})$ for its
non-descendants.
The non-descendants decompose as the disjoint union $\mathrm{ND}_G(k)=\mathrm{anc}_p(k) \cup \mathrm{inc}_p(k)$.
We say that $\nu$ is \emph{Markov with respect to $G$} if it satisfies
the local Markov property
\[
    X_k \indep_\nu X_{\mathrm{ND}_G(k)\setminus\mathrm{Pa}_G(k)}
    \mid X_{\mathrm{Pa}_G(k)},
    \qquad k\in\cI.
\]

\begin{corollary}[EGRT invariance across topological orders]
    \label{cor:dag_ancestral_screening}
    Let $G$ be a DAG on the countable set $\cI$ with
    every $\mathrm{anc}_p(k)$ finite, and suppose $\nu$ is Markov with
    respect to $G$.  Then ancestral screening,
    condition~\textup{(i)} of~\cref{thm:exact_index_order_invariance}, holds
    for the reachability order $\trianglelefteq_p$ of $G$.  Consequently
    every pair of finite-predecessor topological orders of $G$
    produces the same EGRT $\nu \otimes \lambda^{\cI}$-a.e.,
    as well as the same output law under every $\mu \ll \nu$.
    Its $k$-th coordinate is
    \[
        \mathrm{rPIT}^{\preceq_k}_{
            Q_{k \mid \mathrm{Pa}_G(k)}^\nu(x_{\mathrm{Pa}_G(k)},\cdot)}(x_k,v_k)
        \qquad \nu \otimes \lambda^{\cI}\text{-a.e.}
    \]
\end{corollary}

\section{Statistical equivalence across orders}
\label{sec:statistical_equivalence}

\Cref{sec:order_invariance} characterizes when the choice of orders leaves the
EGRT unchanged.
When those conditions fail, the EGRT can produce different output laws under the
alternative for
different orders, leading to differences in the power of a fixed test.
This raises a natural question: does the choice of orders fundamentally
change the \emph{statistical experiment}---that is,
the information available to us about the data law?

We show that the EGRT preserves the original statistical experiment regardless of order.
For every admissible choice of orders, the observation can be
reconstructed from the output almost surely under any law whose
finite-dimensional marginals are dominated by those of the reference.
Consequently, the transformed experiments are equivalent over every
such family, and a test under
one order can be transferred to another with the same power
function.

For reference-dominated laws, we characterize the output likelihood ratio
to explain how the power of a fixed test can depend on the order even
though every $f$-divergence from the reference is preserved.
We also consider the information available in a finite prefix of the output
and show that its divergence from uniformity equals the divergence between the input
marginals on the same coordinates. Different index orders can therefore
expose different amounts of information in their initial segments,
depending on which coordinates they include. As the prefixes grow to
approximate the index set, these divergences increase to the full divergence.

\subsection{Reconstruction and equivalence of experiments}
\label{subsec:equivalence_test_transfer}

Fix the reference law $\nu$ and a nonempty class $\fO$ of admissible
order specifications $o=(\vec{\preceq},\trianglelefteq)$: each $o$
chooses Borel-embeddable total orders on the component spaces and a
finite-predecessor index order.
Write $T_o$ for the EGRT under $o$ and
\begin{equation*}
    \rho_o^\mu
    :=
    (T_o)_{\#}(\mu \otimes \lambda^{\cI}),
\end{equation*}
for its output law under $\mu$, so that
$\rho_o^\nu=\lambda^{\cI}$ by \cref{thm:grt_uniformity}.
Both component and index orders may vary within $\fO$; no
invariance condition from \cref{sec:order_invariance} is assumed.

We call an input law $\mu$ \emph{locally dominated} by $\nu$ if
$\mu_S\ll\nu_S$ for every finite $S\subseteq\cI$. For any
finite-predecessor index order, this is equivalent to requiring
domination on every finite prefix, since every finite set is
contained in a finite prefix. This condition does not depend on
the order. It follows from $\mu\ll\nu$, and is equivalent to it
when $\cI$ is finite, but can be weaker when $\cI$ is infinite.

\begin{lemma}[Measurable reconstruction from an EGRT output]
    \label{lem:egrt_output_reconstruction}
    Fix $o \in \fO$.
    There exists a measurable reconstruction map, depending only on
    $\nu$ and $o$,
    \[
        G_o:[0,1]^{\cI} \to \cX_{\cI}
    \]
    such that, for every $\mu$ locally dominated by $\nu$,
    \begin{equation}
        \label{eq:main_egrt_reconstruction}
        G_o(T_o(X,V))=X
        \qquad
        \mu \otimes \lambda^{\cI}\text{-a.e.},
    \end{equation}
    and hence $(G_o)_{\#}\rho_o^\mu=\mu$.
\end{lemma}

The reconstruction proceeds in index order, recovering each
coordinate from its output and the already reconstructed
predecessors by conditional quantile inversion. The construction
is given in \cref{sec:supp_statistical_equivalence}.

Some restriction on the input family is necessary. For a single
real-valued component with the usual order and reference
$\nu=\Unif[0,1]$, we have $T_o(2,v)=T_o(3,v)=1$ for every $v$.
Thus $\delta_2$ and $\delta_3$ produce the same output law, and no
single reconstruction map can recover the observation under both.
The common output law is nonuniform, so this is consistent with
\cref{thm:grt_converse_uniformization}, which identifies the
reference by uniformity but does not require distinct alternatives
to produce distinct output laws.

Reconstruction yields equivalence of the original and transformed
statistical experiments
\citep{blackwell1953equivalent,torgersen1991comparison}.
An \emph{experiment} is a family of probability measures on a common
measurable space, indexed by the unknown input law $\mu\in\mathcal M$.
Two experiments are called \emph{equivalent} if each can be obtained
from the other via a $\mu$-independent Markov kernel, in which case
they share the same risk functions for all randomized decision rules
(and, in particular, the same power functions for tests).
Here the forward and reconstruction kernels are
\begin{equation*}
    \begin{aligned}
        \mathsf P_o(x,B)
        &:=
        \int_{[0,1]^{\cI}}
            \mathbf 1\!\left\{T_o(x,v)\in B\right\}
            \,\lambda^{\cI}(\mathrm dv),\\
        \mathsf R_o(u,A)
        &:=
        \mathbf 1\!\left\{G_o(u)\in A\right\},
    \end{aligned}
\end{equation*}
for Borel sets $B\subseteq[0,1]^{\cI}$ and $A\subseteq\cX_{\cI}$.
Writing $(\mu\mathsf P_o)(B):=\int\mathsf P_o(x,B)\,\mu(\mathrm dx)$,
we have $\mu\mathsf P_o=\rho_o^\mu$, while reconstruction gives
\begin{equation}
    \label{eq:reverse_kernel_recovers_input_law}
    (\mu\mathsf P_o)\mathsf R_o=\mu
    \qquad
    \text{for every }\mu\text{ locally dominated by }\nu.
\end{equation}
In particular, $(G_o)_{\#}\lambda^{\cI}=\nu$.

\begin{theorem}[Equivalence of the original and transformed experiments]
    \label{thm:exact_testing_equivalence}
    Let $\mathcal M$ be any family of probability measures on
    $(\cX_{\cI},\cB_{\cI})$ such that $\mu$ is locally dominated by $\nu$ for every
    $\mu \in \mathcal M$, and let $o,o' \in \fO$. Then the original
    experiment $(\mu)_{\mu\in\mathcal M}$ and the transformed
    experiments $(\rho_o^\mu)_{\mu\in\mathcal M}$ and
    $(\rho_{o'}^\mu)_{\mu\in\mathcal M}$ are equivalent.
\end{theorem}

A test $\psi$ on the output under $o$, applied after transforming
$X=x$, rejects with conditional probability
\begin{equation*}
    (\mathsf P_o\psi)(x)
    :=
    \int_{[0,1]^{\cI}}\psi(u)\,\mathsf P_o(x,\mathrm du)
    =
    \int_{[0,1]^{\cI}}
        \psi\bigl(T_o(x,v)\bigr)\,\lambda^{\cI}(\mathrm dv).
\end{equation*}
Reconstructing $x$ from the output under $o'$ therefore gives the
following test on that output.

\begin{corollary}[Transfer of tests between orders]
    \label{cor:pairwise_test_transport}
    Fix $o,o'\in\fO$ and a measurable test
    $\psi:[0,1]^{\cI}\to[0,1]$ for the output under $o$.
    Define its transfer to order $o'$ by
    \begin{equation}
        \label{eq:pairwise_test_transport}
        \psi'(u')
        :=(\mathsf P_o\psi)(G_{o'}(u'))
        =\int_{[0,1]^{\cI}}
            \psi\!\left(T_o(G_{o'}(u'),v)\right)
            \,\lambda^{\cI}(\mathrm dv).
    \end{equation}
    Then $\psi'$ is a measurable test with values in $[0,1]$.
    For every $\mu$ locally dominated by $\nu$,
    \begin{equation}
        \label{eq:pairwise_test_transport_power}
        \int\psi'\,\mathrm d\rho_{o'}^\mu
        =
        \int\psi\,\mathrm d\rho_o^\mu.
    \end{equation}
\end{corollary}

Equivalently, reconstruct $x=G_{o'}(u')$, draw fresh independent
auxiliary uniforms, and apply $\psi$ to $T_o(x,V)$. This has
conditional rejection probability $\psi'(u')$, and the transfer depends only on $\nu$, $o$, $o'$, and $\psi$.
The resulting test matches the power of $\psi$ under every law
locally dominated by $\nu$, although its rejection region may
have a different form, which we demonstrate in the following example.

\begin{example}[A tail test with a disconnected transfer]
    \label{ex:order_test_transfer}
    We return to \cref{ex:graph_component_ambiguity}, allowing more
    edges and unequal reference probabilities.
    Consider a single graph-valued component on a fixed vertex set,
    with $m\geq2$ possible edges $e_1,\ldots,e_m$, and write
    $z_j(H):=\mathbf1\{e_j\in H\}$. Under the reference law $\nu$,
    let the edge indicators be independent, with
    $\nu\{z_1=1\}=\alpha\in(0,1)$ and
    $\nu\{z_j=1\}=1/2$ for $j\geq2$.
    Put $N:=2^{m-1}$ and let
    $B(H):=\sum_{j=2}^m2^{m-j}z_j(H)$ be the binary rank of the
    remaining edges. Compare the orders $o,o'$ obtained by sorting
    the graphs by increasing values of
    \[
        h_o(H):=Nz_1(H)+B(H),
        \qquad
        h_{o'}(H):=2B(H)+z_1(H).
    \]
    Both orders extend graph inclusion, but the first prioritizes
    the presence of $e_1$, whereas the second prioritizes every other edge.

    Under $o$, the graphs containing $e_1$ occupy the output interval
    $(1-\alpha,1)$. Thus the upper-tail test
    $\psi(u)=\mathbf1\{u>1-\alpha\}$ has size $\alpha$ and rejects
    exactly when $e_1$ is present, almost surely. Under $o'$, each
    value $b$ of $B$ occupies a block of length $1/N$, and the graphs
    containing $e_1$ occupy the final $\alpha$ fraction of each block.
    The transfer in \cref{eq:pairwise_test_transport} is therefore
    \[
        \psi'(u')
        =\mathbf1\!\left\{
            u'\in\bigcup_{b=0}^{N-1}
            \left(\frac{b+1-\alpha}{N},\frac{b+1}{N}\right)
        \right\},
    \]
    up to a Lebesgue-null set. Under every law $\eta$ on this graph
    space, $\psi$ applied under $o$ and $\psi'$ applied under $o'$
    both reject with probability $\eta\{z_1=1\}$.

    If $N\alpha$ is an integer, take
    $\mu=\nu(\cdot\mid z_1=1)$. The interval $(1-\alpha,1)$
    contains exactly $N\alpha$ complete blocks under $o'$, so
    \[
        \int\psi\,\mathrm d\rho_o^\mu=1,
        \qquad
        \int\psi\,\mathrm d\rho_{o'}^\mu=\alpha,
        \qquad
        \int\psi'\,\mathrm d\rho_{o'}^\mu=1.
    \]
    The same upper-tail test has power equal to its size under
    $o'$, whereas the transferred test retains power one.
\end{example}

\subsection{Likelihood ratios and divergence preservation}
\label{subsec:divergence_preservation}

Fix $o\in\fO$ and $\mu\ll\nu$, and write
$R:=\mathrm d\mu/\mathrm d\nu$ for the input likelihood ratio.
The output likelihood ratio is obtained by composing $R$ with
the reconstruction map.

\begin{proposition}[Output likelihood ratio]
    \label{prop:egrt_output_likelihood_ratio}
    We have $\rho_o^\mu\ll\lambda^{\cI}$ and
    \begin{equation}
        \label{eq:egrt_output_likelihood_ratio}
        \frac{\mathrm d\rho_o^\mu}{\mathrm d\lambda^{\cI}}(u)
        =R(G_o(u))
        \qquad\lambda^{\cI}\text{-a.e.}
    \end{equation}
\end{proposition}

Because $(G_o)_{\#}\lambda^{\cI}=\nu$, the random variable
$R(G_o(U))$ for $U\sim\lambda^{\cI}$ has the same distribution
as $R(X)$ for $X\sim\nu$. This distribution does not depend on
the order, even though the likelihood-ratio function
$u\mapsto R(G_o(u))$ can change.

This also preserves every $f$-divergence
\citep{ali1966general,csiszar1967information}. For probability measures
$P\ll Q$ and a convex function $f:(0,\infty)\to\R$ with $f(1)=0$,
extend $f$ to zero by
$f(0):=\lim_{t\downarrow0}f(t)\in\R\cup\{+\infty\}$ and define
\[
    D_f(P\|Q):=\int f\!\left(\frac{\mathrm dP}{\mathrm dQ}\right)
    \,\mathrm dQ.
\]
An $f$-divergence thus depends only on the distribution of the
likelihood ratio under the reference measure $Q$.

\begin{corollary}[Preservation of divergences]
    \label{cor:egrt_divergence_preservation}
    For every $f$ as above,
    \begin{equation}
        \label{eq:egrt_divergence_preservation}
        D_f\!\left(\rho_o^\mu\middle\|\lambda^{\cI}\right)
        =D_f(\mu\|\nu),
    \end{equation}
    including when both sides are infinite.
\end{corollary}

Applying the data-processing inequality for
$f$-divergences~\citep[Theorem~3.1]{liese2012divergences} to the
forward kernel gives one inequality, and applying it to the
reconstruction kernel gives the reverse, so the divergence is preserved exactly.
Taking $f(t)=|t-1|/2$,
$t\log t$, $(\sqrt t-1)^2$, or $(t-1)^2$ gives preservation of total
variation, Kullback--Leibler divergence, squared Hellinger distance
(without a factor $1/2$), or $\chi^2$ divergence, respectively.

For a fixed test $\psi$, however, its power is
\[
    \int\psi(u)R(G_o(u))\,\lambda^{\cI}(\mathrm du),
\]
which depends on the joint distribution of $\psi$ and $R\circ G_o$
under the uniform law. The power of $\psi$ can therefore vary across orders, even though
the likelihood ratio has the same marginal distribution under
the uniform law.

\begingroup
\theoremstyle{definition}
\newtheorem*{continuedexample}{Example~\ref{ex:order_test_transfer} (continued)}
\begin{continuedexample}
    For the alternative $\mu=\nu(\cdot\mid z_1=1)$, the two output densities are
    \[
        \frac{\mathrm d\rho_o^\mu}{\mathrm d\lambda}(u)
        =\frac{1}{\alpha}\psi(u),
        \qquad
        \frac{\mathrm d\rho_{o'}^\mu}{\mathrm d\lambda}(u)
        =\frac{1}{\alpha}\psi'(u),
    \]
    up to Lebesgue-null sets. Under the uniform law, both densities take
    the value $1/\alpha$ with probability $\alpha$ and $0$ with
    probability $1-\alpha$.
    Thus, for either order,
    \[
        D_f(\rho_o^\mu\|\lambda)
        =D_f(\rho_{o'}^\mu\|\lambda)
        =\alpha f(1/\alpha)+(1-\alpha)f(0)
        =D_f(\mu\|\nu).
    \]
    In particular, both Kullback--Leibler divergences are
    $\log(1/\alpha)$, although one output density is concentrated on an
    upper tail and the other on many separated intervals.
\end{continuedexample}
\endgroup

\subsection{Information in finite prefixes}
\label{subsec:finite_prefix_information}

Fix $o\in\fO$ and an input law $\mu$ locally dominated by $\nu$.
Let $D\subseteq\cI$ be a finite
prefix of the index order in $o$, and write
$\rho_{o,D}^\mu:=(\pi_D)_{\#}\rho_o^\mu$ for the marginal of the
output law on the coordinates in $D$. By projective consistency,
this marginal is the output of the finite-prefix EGRT applied to
$\mu_D$. The divergence identity therefore applies to the input
and output marginals.

\begin{corollary}[Divergences along finite prefixes]
    \label{cor:egrt_finite_prefix_divergence}
    For every finite prefix $D$,
    \begin{equation}
        \label{eq:egrt_finite_prefix_divergence}
        D_f\!\left(\rho_{o,D}^\mu\middle\|\lambda^D\right)
        =D_f(\mu_D\|\nu_D).
    \end{equation}
\end{corollary}

\cref{cor:egrt_finite_prefix_divergence} provides an exact quantification
of the claim in \cref{cor:finite_witness_misspecification}, guaranteeing
not only that every locally dominated
departure from $\nu$ is visible in some finite prefix, but also that the magnitude
of departure matches that in the relevant indices of the original law.
When $\mu\ll\nu$, these divergences also recover the full divergence
as the finite prefixes exhaust the index set.

\begin{corollary}[Recovery of the full divergence]
    \label{cor:egrt_finite_prefix_divergence_limit}
    Suppose in addition that $\mu\ll\nu$, and let $(D_m)_{m\geq1}$ be an increasing sequence of finite prefixes
    with $\bigcup_mD_m=\cI$. Then
    \begin{equation}
        \label{eq:egrt_finite_prefix_divergence_limit}
        D_f\!\left(\rho_{o,D_m}^\mu\middle\|\lambda^{D_m}\right)
        \uparrow D_f(\mu\|\nu)
        \qquad\text{in }[0,\infty].
    \end{equation}
\end{corollary}

\begin{example}[A single-coordinate mean shift under different orders]
    \label{ex:gaussian_prefix_information}
    Take $\cI=\{1,\ldots,d\}$ with $d\geq2$ and the usual order on
    each real-valued component. Let
    \[
        \nu=N(0,I_d),\qquad \mu=N(\delta e_j,I_d),
    \]
    where $e_j$ is the $j$-th standard basis vector and $\delta\neq0$.
    The alternative differs from the reference only through a mean
    shift in coordinate $j$. The coordinates are independent under
    $\nu$, so every index order gives
    $[T_o(x,v)]_k=\Phi(x_k)$ up to null sets, where $\Phi$ is the
    standard normal CDF.

    For an order $o$, write $D_m$ for the set of its first $m$ indices
    and $r_o(j)$ for the position of coordinate $j$.
    For $j \notin D_m$ we have $\mu_{D_m}=\nu_{D_m}$, whereas once
    $j \in D_m$ the marginal likelihood ratio becomes
    $\exp(\delta x_j-\delta^2/2)$, with $\mu$-expected log-likelihood
    ratio $\delta^2/2$.
    Hence, \cref{cor:egrt_finite_prefix_divergence} gives
    \[
        D_{\mathrm{KL}}\!\left(
            \rho_{o,D_m}^\mu\middle\|\lambda^{D_m}\right)
        =D_{\mathrm{KL}}(\mu_{D_m}\|\nu_{D_m})
        =\begin{cases}
            0, & m<r_o(j),\\
            \delta^2/2, & m\geq r_o(j).
        \end{cases}
    \]
    Putting $j$ first makes the full divergence available from the
    first output coordinate alone. Putting it last leaves every
    proper prefix uniform. This is consistent with map invariance,
    since the labeled output does not change. What can change is
    the set $D_m$ of coordinates observed by a test that uses the
    first $m$ outputs.
\end{example}

\section{Discussion}
\label{sec:discussion}

The EGRT provides a common basis for goodness-of-fit testing on
countable products of standard Borel spaces and, for fixed orders,
is unique among locally randomized triangular uniformizing maps
that are lexicographically monotone in the current observation and randomizer.
The choice of order is an important design consideration.
We show conditions under which the EGRT is invariant to this choice,
but they are not guaranteed in practice.

Over reference-dominated families, and more generally families whose
finite-dimensional marginals are dominated by those of the reference,
every order retains the same statistical information available to downstream testing,
and we show how tests can be transferred between orders to
produce equivalent diagnostics. However, transferring a test between orders may make its rejection rule
more complicated or require additional coordinates.

The practical question then is how to choose an order that
makes plausible departures from the reference model easy to
detect with simple diagnostics.
Choosing an order or fitting
the reference law using the observations being tested requires
calibration beyond the fixed-reference, fixed-order results here.
While one can mitigate
sensitivity to any single ordering by using order-invariant test
statistics~\citep{dovern2020order} or pooling evidence across
multiple orders~\citep{pournaderi2026calibrated}, a systematic
framework for selecting or combining orders to target specific
alternatives remains an open problem.

\section*{AI Disclosure}

The authors used GPT-6 Astra to assist with proofreading the prose and checking mathematical results.
The authors take full responsibility for the content of the manuscript, including the correctness of all mathematical arguments.

\bibliographystyle{plainnat}
\bibliography{references}

\begin{thebibliography}{33}
\providecommand{\natexlab}[1]{#1}
\providecommand{\url}[1]{\texttt{#1}}
\expandafter\ifx\csname urlstyle\endcsname\relax
  \providecommand{\doi}[1]{doi: #1}\else
  \providecommand{\doi}{doi: \begingroup \urlstyle{rm}\Url}\fi

\bibitem[Akbari et~al.(2023)Akbari, Ganassali, and
  Kiyavash]{akbari2023learning}
Sina Akbari, Luca Ganassali, and Negar Kiyavash.
\newblock Learning causal graphs via monotone triangular transport maps.
\newblock arXiv preprint arXiv:2305.18210, 2023.
\newblock URL \url{https://arxiv.org/abs/2305.18210}.

\bibitem[Ali and Silvey(1966)]{ali1966general}
S.~M. Ali and S.~D. Silvey.
\newblock A general class of coefficients of divergence of one distribution
  from another.
\newblock \emph{Journal of the Royal Statistical Society: Series B
  (Methodological)}, 28\penalty0 (1):\penalty0 131--142, 1966.
\newblock \doi{10.1111/j.2517-6161.1966.tb00626.x}.

\bibitem[Backhoff et~al.(2017)Backhoff, Beiglb{\"o}ck, Lin, and
  Zalashko]{backhoff2017causal}
Julio Backhoff, Mathias Beiglb{\"o}ck, Yiqing Lin, and Anastasiia Zalashko.
\newblock Causal transport in discrete time and applications.
\newblock \emph{SIAM Journal on Optimization}, 27\penalty0 (4):\penalty0
  2528--2562, 2017.
\newblock \doi{10.1137/16M1080197}.
\newblock URL \url{https://arxiv.org/abs/1606.04062v2}.

\bibitem[Baptista et~al.(2024{\natexlab{a}})Baptista, Hosseini, Kovachki, and
  Marzouk]{baptista2024conditional}
Ricardo Baptista, Bamdad Hosseini, Nikola~B. Kovachki, and Youssef~M. Marzouk.
\newblock Conditional sampling with monotone {GANs}: From generative models to
  likelihood-free inference.
\newblock \emph{SIAM/ASA Journal on Uncertainty Quantification}, 12\penalty0
  (3):\penalty0 868--900, 2024{\natexlab{a}}.

\bibitem[Baptista et~al.(2024{\natexlab{b}})Baptista, Marzouk, and
  Zahm]{baptista2024representation}
Ricardo Baptista, Youssef Marzouk, and Olivier Zahm.
\newblock On the representation and learning of monotone triangular transport
  maps.
\newblock \emph{Foundations of Computational Mathematics}, 24\penalty0
  (6):\penalty0 2063--2108, 2024{\natexlab{b}}.

\bibitem[Beiglb{\"o}ck et~al.(2025)Beiglb{\"o}ck, Pammer, and
  Schrott]{beiglboeck2025denseness}
Mathias Beiglb{\"o}ck, Gudmund Pammer, and Stefan Schrott.
\newblock Denseness of biadapted {Monge} mappings.
\newblock \emph{Annales de l'Institut Henri Poincar{\'e}, Probabilit{\'e}s et
  Statistiques}, 61\penalty0 (1):\penalty0 329--349, 2025.
\newblock \doi{10.1214/23-AIHP1431}.
\newblock URL \url{https://arxiv.org/abs/2210.15554v3}.
\newblock Corrected version: arXiv:2210.15554v3.

\bibitem[Berg(2009)]{berg2009copula}
Daniel Berg.
\newblock Copula goodness-of-fit testing: an overview and power comparison.
\newblock \emph{The European Journal of Finance}, 15\penalty0 (7--8):\penalty0
  675--701, 2009.
\newblock \doi{10.1080/13518470802697428}.

\bibitem[Blackwell(1953)]{blackwell1953equivalent}
David Blackwell.
\newblock Equivalent comparisons of experiments.
\newblock \emph{The Annals of Mathematical Statistics}, 24\penalty0
  (2):\penalty0 265--272, 1953.
\newblock \doi{10.1214/aoms/1177729032}.

\bibitem[Blackwell and Dubins(1975)]{blackwell1975}
David Blackwell and Lester~E. Dubins.
\newblock On existence and non-existence of proper, regular, conditional
  distributions.
\newblock \emph{The Annals of Probability}, 3\penalty0 (5):\penalty0 741--752,
  1975.
\newblock \doi{10.1214/aop/1176996261}.

\bibitem[Bogachev et~al.(2005)Bogachev, Kolesnikov, and
  Medvedev]{bogachev2005triangular}
Vladimir~I. Bogachev, Alexander~V. Kolesnikov, and Kirill~V. Medvedev.
\newblock Triangular transformations of measures.
\newblock \emph{Sbornik: Mathematics}, 196\penalty0 (3):\penalty0 309--335,
  2005.
\newblock \doi{10.1070/SM2005v196n03ABEH000882}.

\bibitem[Brockwell(2007)]{brockwell2007universal}
Anthony~E. Brockwell.
\newblock Universal residuals: A multivariate transformation.
\newblock \emph{Statistics \& Probability Letters}, 77\penalty0 (14):\penalty0
  1473--1478, 2007.
\newblock \doi{10.1016/j.spl.2007.02.008}.

\bibitem[Csisz{\'a}r(1967)]{csiszar1967information}
Imre Csisz{\'a}r.
\newblock Information-type measures of difference of probability distributions
  and indirect observations.
\newblock \emph{Studia Scientiarum Mathematicarum Hungarica}, 2:\penalty0
  299--318, 1967.

\bibitem[Diebold et~al.(1998)Diebold, Gunther, and Tay]{diebold1998}
Francis~X. Diebold, Todd~A. Gunther, and Anthony~S. Tay.
\newblock Evaluating density forecasts with applications to financial risk
  management.
\newblock \emph{International Economic Review}, 39\penalty0 (4):\penalty0
  863--883, 1998.
\newblock \doi{10.2307/2527342}.

\bibitem[Diebold et~al.(1999)Diebold, Hahn, and Tay]{diebold1999multivariate}
Francis~X. Diebold, Jinyong Hahn, and Anthony~S. Tay.
\newblock Multivariate density forecast evaluation and calibration in financial
  risk management: High-frequency returns on foreign exchange.
\newblock \emph{The Review of Economics and Statistics}, 81\penalty0
  (4):\penalty0 661--673, 1999.
\newblock \doi{10.1162/003465399558526}.

\bibitem[Dovern and Manner(2020)]{dovern2020order}
Jonas Dovern and Hans Manner.
\newblock Order-invariant tests for proper calibration of multivariate density
  forecasts.
\newblock \emph{Journal of Applied Econometrics}, 35\penalty0 (4):\penalty0
  440--456, 2020.
\newblock \doi{10.1002/jae.2755}.

\bibitem[Ferguson(1967)]{ferguson1967mathematical}
Thomas~S. Ferguson.
\newblock \emph{Mathematical Statistics: A Decision Theoretic Approach}.
\newblock Academic Press, New York, 1967.

\bibitem[Genest et~al.(2009)Genest, R\'emillard, and
  Beaudoin]{genest2009goodness}
Christian Genest, Bruno R\'emillard, and David Beaudoin.
\newblock Goodness-of-fit tests for copulas: A review and a power study.
\newblock \emph{Insurance: Mathematics and Economics}, 44\penalty0
  (2):\penalty0 199--213, 2009.
\newblock \doi{10.1016/j.insmatheco.2007.10.005}.

\bibitem[Hack et~al.(2022)Hack, Braun, and Gottwald]{hack2022injective}
Pedro Hack, Daniel~A. Braun, and Sebastian Gottwald.
\newblock Representing preorders with injective monotones.
\newblock \emph{Theory and Decision}, 93\penalty0 (4):\penalty0 663--690, 2022.
\newblock \doi{10.1007/s11238-021-09861-w}.

\bibitem[Hack et~al.(2023)Hack, Braun, and Gottwald]{hack2023classification}
Pedro Hack, Daniel~A. Braun, and Sebastian Gottwald.
\newblock The classification of preordered spaces in terms of monotones:
  Complexity and optimization.
\newblock \emph{Theory and Decision}, 94:\penalty0 693--720, 2023.
\newblock \doi{10.1007/s11238-022-09904-w}.

\bibitem[Justel et~al.(1997)Justel, Pe\~na, and Zamar]{justel1997multivariate}
Ana Justel, Daniel Pe\~na, and Rub\'en Zamar.
\newblock A multivariate {Kolmogorov--Smirnov} test of goodness of fit.
\newblock \emph{Statistics \& Probability Letters}, 35\penalty0 (3):\penalty0
  251--259, 1997.
\newblock \doi{10.1016/S0167-7152(97)00020-5}.

\bibitem[Kallenberg(2021)]{kallenberg2021foundations}
Olav Kallenberg.
\newblock \emph{Foundations of Modern Probability}.
\newblock Springer, third edition, 2021.
\newblock \doi{10.1007/978-3-030-61871-1}.

\bibitem[Kechris(1995)]{kechris1995classical}
Alexander~S. Kechris.
\newblock \emph{Classical Descriptive Set Theory}, volume 156 of \emph{Graduate
  Texts in Mathematics}.
\newblock Springer, 1995.
\newblock \doi{10.1007/978-1-4612-4190-4}.

\bibitem[Knothe(1957)]{knothe1957}
Herbert Knothe.
\newblock Contributions to the theory of convex bodies.
\newblock \emph{Michigan Mathematical Journal}, 4\penalty0 (1):\penalty0
  39--52, 1957.
\newblock \doi{10.1307/mmj/1028990175}.

\bibitem[Lauritzen(1996)]{lauritzen1996graphical}
Steffen~L. Lauritzen.
\newblock \emph{Graphical Models}, volume~17 of \emph{Oxford Statistical
  Science Series}.
\newblock Oxford University Press, 1996.

\bibitem[Liese(2012)]{liese2012divergences}
Friedrich Liese.
\newblock $\phi$-divergences, sufficiency, {Bayes} sufficiency, and deficiency.
\newblock \emph{Kybernetika}, 48\penalty0 (4):\penalty0 690--713, 2012.
\newblock URL \url{https://www.kybernetika.cz/content/2012/4/690/paper.pdf}.

\bibitem[Pournaderi(2026)]{pournaderi2026calibrated}
Mehrdad Pournaderi.
\newblock Calibrated order-randomized {R}osenblatt tests.
\newblock arXiv preprint arXiv:2609.37477, 2026.
\newblock URL \url{https://arxiv.org/abs/2609.37477}.

\bibitem[Rosenblatt(1952)]{rosenblatt1952remarks}
Murray Rosenblatt.
\newblock Remarks on a multivariate transformation.
\newblock \emph{The Annals of Mathematical Statistics}, 23\penalty0
  (3):\penalty0 470--472, 1952.
\newblock \doi{10.1214/aoms/1177729394}.

\bibitem[R{\"u}schendorf(1981)]{ruschendorf1981stochastically}
Ludger R{\"u}schendorf.
\newblock Stochastically ordered distributions and monotonicity of the
  {OC}-function of sequential probability ratio tests.
\newblock \emph{Mathematische Operationsforschung und Statistik. Series
  Statistics}, 12\penalty0 (3):\penalty0 327--338, 1981.
\newblock \doi{10.1080/02331888108801593}.

\bibitem[R\"uschendorf(2009)]{ruschendorf2009distributional}
Ludger R\"uschendorf.
\newblock On the distributional transform, {Sklar's} theorem, and the empirical
  copula process.
\newblock \emph{Journal of Statistical Planning and Inference}, 139\penalty0
  (11):\penalty0 3921--3927, 2009.
\newblock \doi{10.1016/j.jspi.2009.05.030}.

\bibitem[R\"uschendorf and de~Valk(1993)]{ruschendorf1993regression}
Ludger R\"uschendorf and Vincent de~Valk.
\newblock On regression representations of stochastic processes.
\newblock \emph{Stochastic Processes and their Applications}, 46\penalty0
  (2):\penalty0 183--198, 1993.
\newblock \doi{10.1016/0304-4149(93)90001-K}.

\bibitem[Skorokhod(1976)]{skorokhod1976representation}
A.~V. Skorokhod.
\newblock On a representation of random variables.
\newblock \emph{Theory of Probability and Its Applications}, 21\penalty0
  (3):\penalty0 628--632, 1976.
\newblock \doi{10.1137/1121076}.

\bibitem[Szpilrajn(1930)]{szpilrajn1930extension}
Edward Szpilrajn.
\newblock Sur l'extension de l'ordre partiel.
\newblock \emph{Fundamenta Mathematicae}, 16\penalty0 (1):\penalty0 386--389,
  1930.
\newblock \doi{10.4064/fm-16-1-386-389}.

\bibitem[Torgersen(1991)]{torgersen1991comparison}
Erik Torgersen.
\newblock \emph{Comparison of Statistical Experiments}, volume~36 of
  \emph{Encyclopedia of Mathematics and its Applications}.
\newblock Cambridge University Press, 1991.
\newblock \doi{10.1017/CBO9780511666353}.

\end{thebibliography}

%%%%%%%%%%%%%%%%%%%%%%%%%%%%%%%%%%%%%%%%%%%%%%%%%%%%%%%%%%%%%%%%%%%%%%%%%%%%%%%%
% Supplementary material

\clearpage
\setcounter{page}{1}
\setcounter{section}{0}
\setcounter{table}{0}
\setcounter{figure}{0}
\renewcommand{\theHsection}{SIsection.\arabic{section}}
\renewcommand{\theHtable}{SItable.\arabic{section}.\arabic{table}}
\renewcommand{\theHfigure}{SIfigure.\arabic{section}.\arabic{figure}}
\renewcommand{\thepage}{S\arabic{page}}
\renewcommand{\thesection}{S\arabic{section}}
\renewcommand{\thetable}{S\arabic{section}.\arabic{table}}
\renewcommand{\thefigure}{S\arabic{section}.\arabic{figure}}

\begin{center}
{\LARGE Supplementary material}
\end{center}

\label{supplement}

\section{Notation}
\label{sec:supp_notation}

The table below collects the notation used in the main text and supplement.
\bigskip

{%
\small
\begin{longtable}{@{}p{4.15cm} p{9.05cm}@{}}
\toprule
\textbf{Symbol} & \textbf{Meaning} \\
\midrule
\endfirsthead
\toprule
\textbf{Symbol} & \textbf{Meaning} \\
\midrule
\endhead
\bottomrule
\endlastfoot

\multicolumn{2}{@{}l}{\textit{Spaces, measures, and kernels}} \\[2pt]
$(\cX_k,\cB_k)$ & Standard Borel component space at coordinate $k$. \\
$(\cX_S,\cB_S)$ & Product space $\prod_{k \in S}\cX_k$ with product $\sigma$-algebra, for $S \subseteq \cI$; $(\cX_{\cI},\cB_{\cI})$ denotes the full countable product. \\
$x_S,\,X_S$ & Restriction of a point $x_{\cI}$, or of the random element $X_{\cI}$, to coordinates in $S$. \\
$\nu,\,\mu$ & Reference and input probability measures on $\cX_{\cI}$, respectively. \\
$\lambda,\,\lambda^S$ & Lebesgue measure on $[0,1]$ and its product on $[0,1]^S$; $\lambda^{\cI}$ is the auxiliary-uniform law. \\
$Q_k^\nu(x_{\triangleleft k},\cdot)$ & Fixed version of the conditional law of $X_k$ given its strict predecessors under the current finite-predecessor index order. \\
$Q_{k \mid S}^\nu(x_S,\cdot)$ & Fixed version of the conditional law of $X_k$ given $X_S$, for $S \subseteq \cI \setminus \{k\}$. \\
$\cH_k^\nu$ & Family $\{Q_k^\nu(h,\cdot):h\in\cX_{\triangleleft k}\}$ for the fixed kernel version, including histories in $\nu_{\triangleleft k}$-null sets. \\
$\nu_S,\,\mu_S$ & Marginals of $\nu$ and $\mu$ on coordinates $S \subseteq \cI$. \\
$\Law_\mu(Z)$ & Law of $Z$ when $X_{\cI} \sim \mu$, with independent auxiliary uniforms understood. \\
\addlinespace

\multicolumn{2}{@{}l}{\textit{Component-space orders}} \\[2pt]
$\preceq_k,\, \prec_k$ & Borel-embeddable total order on $\cX_k$ and its strict part. \\
$h_{\preceq_k}$ & Borel order embedding of $(\cX_k, \preceq_k)$ into $(\R, \leq)$. \\
$\preceq_{p,k}$ & Partial order on $\cX_k$. \\
$\vec{\preceq}=(\preceq_k)_{k \in \cI}$ & Sequence of component-space total orders. \\
$\cE(\preceq_{p,k})$ & Class of Borel-embeddable linear extensions of $\preceq_{p,k}$. \\
$I_{p,k}(x)$ & Set of elements of $\cX_k$ incomparable with $x$ under $\preceq_{p,k}$. \\
$\Delta_k(x_{\trianglelefteq k})$ & Conditional ambiguity $Q_k^\nu(x_{\triangleleft k},I_{p,k}(x_k))$. \\
$(f_n)_{n\geq1}$ & Countable Borel monotone functions: point-separating in \cref{thm:admissibility_characterization}; representing the partial order through \cref{eq:countable_order_representation} in \cref{thm:countable_representation_cdf_separating}. \\
\addlinespace

\multicolumn{2}{@{}l}{\textit{Index-set orders and finite input blocks}} \\[2pt]
$\trianglelefteq,\, \triangleleft$ & Total index order with finite strict predecessor sets, and its strict part. \\
$P_{\triangleleft}(k)$ & Strict predecessor set $\{j \in \cI:j \triangleleft k\}$. \\
$x_{\triangleleft k},\,x_{\trianglelefteq k}$ & Coordinates strictly preceding $k$, and those preceding or equal to $k$, under the current index order. \\
$\trianglelefteq_p$ & Partial order on the index set $\cI$. \\
$\mathcal L(\trianglelefteq_p)$ & Finite-predecessor total-order linear extensions of $\trianglelefteq_p$. \\
$\mathrm{anc}_p(k)$ & Strict ancestors of $k$ under $\trianglelefteq_p$. \\
$\mathrm{inc}_p(k)$ & Indices incomparable with $k$ under $\trianglelefteq_p$. \\
$\mathrm{Pa}_G(k),\,\mathrm{ND}_G(k)$ & Parents and non-descendants of $k$ in a DAG $G$ on $\cI$. \\
\addlinespace

\multicolumn{2}{@{}l}{\textit{Transforms}} \\[2pt]
$F_\eta^\preceq(x),\,F_\eta^\preceq(x^-)$ & CDF and strict-left CDF of $\eta$ under $\preceq$. \\
$\mathrm{rPIT}_\eta^\preceq(x,v)$ & Randomized probability integral transform $F_\eta^\preceq(x^-)+v\eta(\{x\})$. \\
$T_{\cI,\nu}^{\vec \preceq, \trianglelefteq}$ & EGRT constructed from reference law $\nu$, component orders $\vec \preceq$, and index order $\trianglelefteq$. \\
$U^{\vec \preceq, \trianglelefteq}$ & EGRT output $T_{\cI,\nu}^{\vec \preceq, \trianglelefteq}(X_{\cI},V_{\cI})$. \\
$T_{\cI,\nu}^{\vec \preceq},\,U^{\vec \preceq}$ & Notation when the index order is fixed and dependence on component orders is emphasized. \\
$T_{\cI,\nu}^{\trianglelefteq},\,U^{\trianglelefteq}$ & Notation when component orders are fixed and dependence on the index order is emphasized. \\
\textit{locally randomized triangular} & The $k$-th output may depend on $X_{\trianglelefteq k}$ and the current $V_k$, but not on earlier auxiliary variables. \\
\addlinespace

\multicolumn{2}{@{}l}{\textit{Reconstruction and statistical experiments}} \\[2pt]
$\fO,\,T_o$ & Class of admissible order specifications and the EGRT under $o$. \\
$\mathsf P_o,\,\rho_o^\mu$ & Forward kernel obtained by averaging over auxiliary uniforms, and output law $\mu\mathsf P_o$; $\mathsf P_o\psi$ is the conditional rejection probability of a transformed-space test $\psi$. \\
$G_o,\,\mathsf R_o$ & Measurable reconstruction map and its deterministic Markov kernel. \\
$\mathcal M$ & Family of input laws locally dominated by $\nu$: $\mu_S\ll\nu_S$ for every finite $S\subseteq\cI$. \\
$R$ & Likelihood ratio $\mathrm d\mu/\mathrm d\nu$ for $\mu\ll\nu$. \\
$D_f(P\|Q)$ & $f$-divergence $\int f(\mathrm dP/\mathrm dQ)\,\mathrm dQ$ for $P\ll Q$, with the generator specified in \cref{subsec:divergence_preservation}. \\
$D,\,\rho_{o,D}^\mu$ & Finite downward-closed prefix and the corresponding EGRT output marginal $(\pi_D)_{\#}\rho_o^\mu$; $(D_m)_{m\geq1}$ denotes an increasing sequence of such prefixes. \\
\end{longtable}
}

\setcounter{section}{1}
\section{Proofs for Section~2}
\label{sec:supp_egrt}

\paragraph{Proof of~\cref{lem:disintegration_standard_borel}.}

\begin{proof}
    Applying the standard disintegration theorem
    \citep[Theorem~8.5]{kallenberg2021foundations} to the pair of
    standard Borel random elements $(\pi_{\triangleleft k},\pi_k)$ under
    $\nu$ yields the stated kernel and its
    $\nu_{\triangleleft k}$-a.e.\ uniqueness. When $P_{\triangleleft}(k)$ is empty,
    the conditioning space is a singleton and the kernel reduces to the marginal
    $\nu_{\{k\}}$. Because $(\pi_{\triangleleft k},\pi_k)$ depends on the index order
    only through $P_{\triangleleft}(k)$, the same holds for the kernel up to $\nu_{\triangleleft k}$-null sets.
\end{proof}

\begin{lemma}[Joint measurability of the conditional rPIT]
    \label{lem:kernel_valued_rpit_measurability}
    Let $(S,\cS)$ be a measurable space, $(\cX,\cB)$ a standard Borel space, $K:S \times \cB \to [0,1]$ a probability kernel, and $\preceq$ a Borel-embeddable total order on $\cX$ with Borel order embedding $h:\cX \to \R$.
    Then the map
    \[
        (s,x,v)\longmapsto\mathrm{rPIT}^\preceq_{K(s,\cdot)}(x,v)
        =K(s,\{y \in \cX:y \prec x\})+vK(s,\{x\})
    \]
    from $S \times \cX \times [0,1]$ to $[0,1]$ is $\cS \otimes \cB \otimes \cB([0,1])$-measurable.
\end{lemma}

\begin{proof}
    Let
    \[
        \mathcal M
        :=\left\{A \in \cB \otimes \cB:
        (s,x) \mapsto K(s,A_x)\text{ is }\cS \otimes \cB\text{-measurable}\right\},
    \]
    where $A_x:=\{y \in \cX:(x,y) \in A\}$.
    For any measurable rectangle $A=B_1 \times B_2$, we have
    $K(s,A_x)=\mathbf 1_{B_1}(x)K(s,B_2)$, so $\mathcal M$ contains all measurable rectangles.
    Moreover, $\mathcal M$ is closed under complements since $(A^c)_x=(A_x)^c$ and $K(s,(A^c)_x)=1-K(s,A_x)$, and it is closed under countable disjoint unions because for disjoint $(A_n)_{n \geq 1} \subseteq \mathcal M$, $K(s,(\cup_n A_n)_x)=\sum_n K(s,(A_n)_x)$ is a pointwise limit of measurable partial sums. By the $\pi$--$\lambda$ theorem, $\mathcal M=\cB \otimes \cB$.

    Now consider the sets
    \[
        L:=\{(x,y) \in \cX^2:y \prec x\},
        \qquad
        D:=\{(x,y) \in \cX^2:y=x\}.
    \]
    Because $h$ is Borel, $L=(h \times h)^{-1}\{(a,b) \in \R^2:b<a\}$ is Borel in $\cX^2$, and the diagonal $D$ is Borel since $\cX$ is standard Borel.
    Applying the property defining $\mathcal M$ to $L$ and $D$, we see that both
    \[
        (s,x) \mapsto K(s,L_x)=K(s,\{y:y \prec x\})
        \quad\text{and}\quad
        (s,x) \mapsto K(s,D_x)=K(s,\{x\})
    \]
    are $\cS \otimes \cB$-measurable, from which joint measurability of $(s,x,v)\mapsto \mathrm{rPIT}^\preceq_{K(s,\cdot)}(x,v)$ is immediate.
\end{proof}

\begin{lemma}[Measurable reconstruction from the conditional rPIT]
    \label{lem:conditional_rpit_reconstruction}
    Let $(S,\cS)$ and $(\cX,\cB)$ be standard Borel spaces, $K:S \times \cB \to [0,1]$ a probability kernel, and $\preceq$ a Borel-embeddable total order on $\cX$ with Borel order embedding $h:\cX \to \R$.
    There exists a measurable map $G_K:S \times [0,1] \to \cX$ such that, for every probability measure $\rho$ on $S$, if
    \[
        (S_0,X,V) \sim \rho(ds)K(s,dx)\lambda(dv),
    \]
    then
    \[
        G_K\bigl(S_0,\mathrm{rPIT}^\preceq_{K(S_0,\cdot)}(X,V)\bigr)=X
        \qquad\text{a.s.}
    \]
    Moreover, for each $s \in S$, $(G_K(s,\cdot))_{\#}\lambda=K(s,\cdot)$.
\end{lemma}

\begin{proof}
    By the Lusin--Souslin theorem
    \citep[Corollary~15.2]{kechris1995classical}, the image $Y:=h(\cX)$ is Borel in $\R$ and $h^{-1}:Y \to \cX$ is Borel.
    Let $K^h(s,\cdot):=h_{\#}K(s,\cdot)$ denote the pushforward kernel on $\R$, and define its cumulative distribution function $F_s(t):=K^h(s,(-\infty,t])$, which is measurable in $s$ for each fixed $t$ and right-continuous in $t$ for each fixed $s$.
    Setting
    \[
        H(s,u):=
        \begin{cases}
            \inf\{t \in \Q:F_s(t) \geq u\}, & 0<u<1,\\
            0, & u \in \{0,1\},
        \end{cases}
    \]
    gives a jointly measurable generalized inverse on $S\times[0,1]$ since the infimum is taken over a countable set.
    Fixing an arbitrary point $x_\star \in \cX$, define
    \[
        G_K(s,u):=
        \begin{cases}
            h^{-1}(H(s,u)), & H(s,u) \in Y,\\
            x_\star, & H(s,u)\notin Y.
        \end{cases}
    \]
    For each $s \in S$ and $U \sim \lambda$, we have $H(s,U) \sim K^h(s,\cdot)$.
    Because $K^h(s,Y)=1$, the event $\{H(s,U)\notin Y\}$ is $\lambda$-null, and thus $(G_K(s,\cdot))_{\#}\lambda=K(s,\cdot)$ for all $s$.

    To verify the almost-sure reconstruction property, note that under the embedding $h$,
    \[
        \mathrm{rPIT}^\preceq_{K(s,\cdot)}(x,v)
        =
        F_s(h(x)^-)+vK^h(s,\{h(x)\}).
    \]
    By the standard quantile identity for the real-valued randomized PIT
    \citep[Proposition~2.1]{ruschendorf2009distributional},
    \[
        H\bigl(s,\mathrm{rPIT}^\preceq_{K(s,\cdot)}(X,V)\bigr)=h(X)
        \qquad\text{a.s.}
    \]
    under $K(s,dx)\lambda(dv)$ for each $s\in S$.
    Applying $h^{-1}$ and integrating with respect to $\rho(ds)$ yields the stated identity. Note that $G_K$ is constructed purely from $K$ and $h$, independently of $\rho$.
\end{proof}

\begin{lemma}[Uniqueness of the rPIT under lexicographic monotonicity]
    \label{lem:one_dimensional_monotone_uniqueness}
    Let $(\cX,\cB)$ be a standard Borel space with Borel-embeddable total order $\preceq$, and let $\eta$ be a probability measure on $\cX$.
    Suppose $S:\cX \times [0,1] \to [0,1]$ is measurable, nondecreasing for the lexicographic order on $\cX \times [0,1]$, and satisfies
    \[
        S_{\#}(\eta \otimes \lambda)=\lambda.
    \]
    Then
    \[
        S(x,v)=\mathrm{rPIT}_\eta^\preceq(x,v)
        \qquad
        (\eta \otimes \lambda)\text{-a.e.}
    \]
\end{lemma}

\begin{proof}
    Via a Borel order embedding $h:\cX \to \R$, we may push $\eta$ and $S$ forward to $h(\cX) \subseteq \R$ and thus assume without loss of generality that $\cX \subseteq \R$ with the standard order.
    Write $F(t):=\eta((-\infty,t])$, $F(t^-):=\eta((-\infty,t))$, and $\rho:=\eta \otimes \lambda$.

    Fix $a \in \R$ and put $L_a:=(-\infty,a] \times [0,1]$, so that $\alpha_a:=\rho(L_a)=F(a)$.
    Since every point in $L_a$ lexicographically precedes every point in $L_a^c$, monotonicity of $S$ implies that whenever $\rho(L_a)>0$ and $\rho(L_a^c)>0$,
    \[
        s_a:=\operatorname*{ess\,sup}_{L_a} S
        \leq
        r_a:=\operatorname*{ess\,inf}_{L_a^c} S.
    \]
    We claim that
    \[
        S \leq \alpha_a\quad \rho\text{-a.e.\ on }L_a,
        \qquad
        S \geq \alpha_a\quad \rho\text{-a.e.\ on }L_a^c.
    \]
    (If either $L_a$ or $L_a^c$ is $\rho$-null, the corresponding bound is vacuous and the other is trivial since $S\in[0,1]$.)
    Indeed, if $s_a>\alpha_a$, we can pick $t \in (\alpha_a, s_a)$ with $t \leq r_a$; then $\{S \leq t\} \subseteq L_a$ up to a $\rho$-null set, which gives $\rho\{S \leq t\} \leq \alpha_a < t$, contradicting $S_{\#}\rho=\lambda$.
    Symmetrically, if $r_a<\alpha_a$, picking $t \in (r_a, \alpha_a)$ with $t \geq s_a$ implies that $L_a \subseteq \{S \leq t\}$ up to a null set while $\rho(\{S \leq t\} \cap L_a^c)>0$, so $\rho\{S \leq t\} > \alpha_a > t$, again a contradiction.

    Intersecting these full-measure events over all rational $a \in \Q$, we obtain
    \[
        F(x^-)
        \leq
        S(x,v)
        \leq
        F(x)
        \qquad\text{for }\rho\text{-a.e.\ }(x,v).
    \]
    In particular, $S(x,v)=F(x)=F(x^-)=\mathrm{rPIT}_\eta(x,v)$ for $\rho$-a.e.\ $(x,v)$ such that $\eta(\{x\})=0$.

    It remains to handle the atoms of $\eta$, of which there are at most countably many.
    Fix an atom $x$ and set $a:=F(x^-)$, $m:=\eta(\{x\})>0$, and $b:=F(x)=a+m$.
    By the bounds above, $S(x,v) \in [a,b]$ for $\lambda$-a.e.\ $v$.
    For any $t \in [a,b)$, lexicographic monotonicity yields
    \[
        \{S \leq t\}
        =
        ((-\infty,x) \times [0,1])
        \cup
        (\{x\} \times \{v:S(x,v) \leq t\})
    \]
    up to a $\rho$-null set, and therefore
    \[
        t
        =
        \rho\{S \leq t\}
        =
        a+m\lambda\{v:S(x,v) \leq t\}.
    \]
    That is, $v \mapsto S(x,v)$ pushes $\lambda$ forward to $\mathrm{Unif}(a,b)$.
    Since $v \mapsto S(x,v)$ is nondecreasing, we must have $S(x,v)=a+mv=\mathrm{rPIT}_\eta(x,v)$ for $\lambda$-a.e.\ $v$.
\end{proof}

\paragraph{Proof of~\cref{thm:randomized_PIT_uniform}.}

\begin{proof}
    Let $h:\cX \to \R$ be a Borel order embedding for $\preceq$, and write $\mu^h:=h_{\#}\mu$ and $\nu^h:=h_{\#}\nu$.
    Because
    \[
        \{y \in \cX:y \preceq x\}=h^{-1}((-\infty,h(x)])
        \quad\text{and}\quad
        \{y \in \cX:y \prec x\}=h^{-1}((-\infty,h(x))),
    \]
    we have $\mathrm{rPIT}_\nu^\preceq(x,v) = \mathrm{rPIT}_{\nu^h}(h(x),v)$, where the right-hand side is the standard real-line randomized PIT.
    Under $(X,V) \sim \mu \otimes \lambda$, we have $(h(X),V) \sim \mu^h \otimes \lambda$, so the condition $(\mathrm{rPIT}_\nu^\preceq)_{\#}(\mu \otimes \lambda)=\lambda$ is equivalent to
    \[
        (\mathrm{rPIT}_{\nu^h})_{\#}(\mu^h \otimes \lambda)=\lambda.
    \]
    Moreover, since $h$ is a Borel isomorphism onto its image, $\mu=\nu$ if and only if $\mu^h=\nu^h$.
    When $\mu^h=\nu^h$, the displayed identity holds by standard uniformity of the real-valued randomized PIT~\citep[Proposition~2.1]{ruschendorf2009distributional}.

    Conversely, suppose $(\mathrm{rPIT}_{\nu^h})_{\#}(\mu^h \otimes \lambda)=\lambda$, and let $F_\mu, F_\nu$ denote the CDFs of $\mu^h, \nu^h$.
    Fix $t \in \R$, set $u_0:=F_\nu(t)$, and define
    \[
        A:=(-\infty,t] \times [0,1],
        \qquad
        B:=\{(y,v):\mathrm{rPIT}_{\nu^h}(y,v) \leq u_0\}.
    \]
    Clearly $A \subseteq B$, since $y \leq t$ implies $\mathrm{rPIT}_{\nu^h}(y,v) \leq F_\nu(t)=u_0$.
    On the other hand, if $(y,v) \in B \setminus A$, then $y>t$ and hence $F_\nu(y^-) \geq u_0$, which forces $F_\nu(y^-)=u_0$ and $v\nu^h(\{y\})=0$.
    Among such points $y>t$, those with $\nu^h(\{y\})>0$ require $v=0$ (a $\lambda$-null slice), whereas all $y>t$ with $\nu^h(\{y\})=0$ and $F_\nu(y^-)=u_0$ satisfy $\mathrm{rPIT}_{\nu^h}(y,v)=u_0$ for all $v$.
    Because the pushforward under $\mu^h\otimes\lambda$ is $\lambda$, which is atomless and thus assigns mass zero to $\{u_0\}$, $\mu^h$ must assign mass zero to $\{y>t: \nu^h(\{y\})=0,\, F_\nu(y^-)=u_0\}$.
    Consequently, $(\mu^h \otimes \lambda)(B \setminus A)=0$, and thus
    \[
        F_\mu(t)
        =
        (\mu^h \otimes \lambda)(A)
        =
        (\mu^h \otimes \lambda)(B)
        =
        \lambda([0,u_0])
        =
        F_\nu(t)
    \]
    for all $t\in\R$, proving $\mu^h=\nu^h$.
\end{proof}

For a finite $\trianglelefteq$-downward closed set $D \subseteq \cI$, we write
\[
    T_{D,\nu_D}^{\vec{\preceq}, \trianglelefteq}:
    \cX_D \times [0,1]^D \to [0,1]^D
\]
for the EGRT of $\nu_D$ under the restricted orders, defined using the fixed kernels $Q_k^\nu$ from the main text (which remain valid disintegration kernels for $\nu_D$ since $P_{\triangleleft}(k) \subseteq D$ for all $k \in D$). For $D \subseteq E \subseteq \cI$, let $\pi_D^\cX:\cX_E \to \cX_D$ and $\pi_D^{[0,1]}:[0,1]^E \to [0,1]^D$ denote the coordinate projections.

\paragraph{Proof of~\cref{thm:grt_uniformity}.}

\begin{proof}
    Following the sequential uniformization argument of \citet[Theorem~3]{ruschendorf1993regression}, enumerate $\cI$ in increasing $\trianglelefteq$-order as $k_1,\ldots,k_n$ (or $k_1,k_2,\ldots$ if $\cI$ is infinite) and set $D_m:=\{k_1,\ldots,k_m\}$.
    We show by induction on $m$ that
    \[
        \bigl(T_{D_m,\nu_{D_m}}^{\vec{\preceq}, \trianglelefteq}\bigr)_{\#}(\nu_{D_m} \otimes \lambda^{D_m})
        =
        \lambda^{D_m}.
    \]
    The base case $m=1$ is \cref{thm:randomized_PIT_uniform} applied to $\nu_{\{k_1\}}$.
    For $m \geq 2$, disintegrate $\nu_{D_m}$ as
    \[
        \nu_{D_m}(dx_{D_m})
        =
        \nu_{D_{m-1}}(dx_{D_{m-1}})
        Q_{k_m}^\nu(x_{D_{m-1}},dx_{k_m}).
    \]
    For any measurable rectangle $A \times B \subseteq [0,1]^{D_{m-1}} \times [0,1]$, Fubini's theorem, \cref{thm:randomized_PIT_uniform} applied conditionally given $x_{D_{m-1}}$, and the induction hypothesis give
    \[
    \begin{aligned}
        & (\nu_{D_m} \otimes \lambda^{D_m})
        \left\{T_{D_m,\nu_{D_m}}^{\vec{\preceq}, \trianglelefteq} \in A \times B\right\} \\
        &\qquad =
        \int \mathbf{1}_A\!
        \left(T_{D_{m-1},\nu_{D_{m-1}}}^{\vec{\preceq}, \trianglelefteq}(x_{D_{m-1}},v_{D_{m-1}})\right)
        \lambda(B)\,
        \nu_{D_{m-1}}(dx_{D_{m-1}})\lambda^{D_{m-1}}(dv_{D_{m-1}}) \\
        &\qquad =
        \lambda^{D_{m-1}}(A)\lambda(B).
    \end{aligned}
    \]
    By a monotone-class argument, the pushforward on $D_m$ is $\lambda^{D_m}$.

    This completes the proof when $\cI$ is finite.
    When $\cI$ is countably infinite, any finite subset $S \subseteq \cI$ is contained in
    $D(S):=\bigcup_{k \in S}\{j \in \cI:j \trianglelefteq k\}$, which is finite and downward closed by the finite-predecessor assumption.
    Uniformity on $D(S)$ implies uniformity on $S$ by marginalization, and since $S$ was arbitrary, the full pushforward law is $\lambda^{\cI}$.
\end{proof}

\begin{proposition}[Pointwise projective consistency with inherited kernels]
    \label{prop:grt_projective_consistency}
    Let $D \subseteq E \subseteq \cI$ be finite $\trianglelefteq$-downward closed sets. Then, for every $(x_E,v_E) \in \cX_E \times [0,1]^E$,
    \[
        \pi_D^{[0,1]}\!\left(
            T_{E,\nu_E}^{\vec{\preceq}, \trianglelefteq}(x_E,v_E)
        \right)
        =
        T_{D,\nu_D}^{\vec{\preceq}, \trianglelefteq}
        \bigl(\pi_D^{\cX}x_E,\pi_D^{[0,1]}v_E\bigr).
    \]
    Consequently, for every finite downward-closed $D \subseteq \cI$ and every $(x_{\cI},v_{\cI}) \in \cX_{\cI} \times [0,1]^{\cI}$,
    \[
        \pi_D^{[0,1]}\!\left(
            T_{\cI,\nu}^{\vec{\preceq}, \trianglelefteq}(x_{\cI},v_{\cI})
        \right)
        =
        T_{D,\nu_D}^{\vec{\preceq}, \trianglelefteq}
        \bigl(\pi_D^{\cX}x_{\cI},\pi_D^{[0,1]}v_{\cI}\bigr).
    \]
\end{proposition}

\paragraph{Proof of~\cref{prop:grt_projective_consistency}.}
\begin{proof}
Because $D$ is downward closed, $P_{\triangleleft}(k) \subseteq D$ for every $k \in D$, so $k$ has the exact same predecessor set in $D$, $E$, and $\cI$. Since all three transforms are defined using the same inherited kernel $Q_k^\nu$ and component order $\preceq_k$ at coordinate $k$, their $k$-th coordinates agree pointwise for every $k \in D$.
\end{proof}

\paragraph{Proof of~\cref{thm:grt_converse_uniformization}.}

\begin{proof}
    The forward direction is \cref{thm:grt_uniformity}.
    For the converse, fix a finite nonempty $\trianglelefteq$-downward closed set $D$, enumerated as $k_1 \triangleleft \cdots \triangleleft k_n$, and put $D_m:=\set{k_1,\ldots,k_m}$.
    Assuming
    \[
        \bigl(T_{D,\nu_D}^{\vec{\preceq}, \trianglelefteq}\bigr)_{\#}
        (\mu_D \otimes \lambda^D)=\lambda^D,
    \]
    we show that $\mu_{D_m}=\nu_{D_m}$ by induction on $m$.
    The case $m=1$ follows directly from \cref{thm:randomized_PIT_uniform}.

    Suppose $m \geq 2$ and $\mu_{D_{m-1}}=\nu_{D_{m-1}}$.
    Let $R_m(x_{D_{m-1}},\cdot)$ be a regular conditional distribution of $X_{k_m}$ given $X_{D_{m-1}}$ under $\mu_{D_m}$.
    Applying \cref{lem:conditional_rpit_reconstruction} recursively to the first $m-1$ coordinates yields a measurable map
    $S_{m-1}:[0,1]^{D_{m-1}} \to \cX_{D_{m-1}}$ such that $(S_{m-1})_{\#}\lambda^{D_{m-1}}=\nu_{D_{m-1}}$ and
    $S_{m-1}(U_{D_{m-1}})=X_{D_{m-1}}$ a.s.\ under $\nu_{D_{m-1}} \otimes \lambda^{D_{m-1}}$---and hence also under $\mu_{D_{m-1}} \otimes \lambda^{D_{m-1}}$, since $\mu_{D_{m-1}}=\nu_{D_{m-1}}$ by the induction hypothesis (with no domination assumption required).

    Because $X_{k_m} \indep V_{D_{m-1}} \mid X_{D_{m-1}}$, conditioning first on $(X_{D_{m-1}},V_{D_{m-1}})$ shows that $R_m(S_{m-1}(u),\cdot)$ is a version of the conditional law of $X_{k_m}$ given $U_{D_{m-1}}=u$.
    By hypothesis, $U_{k_m}=\mathrm{rPIT}^{\preceq_{k_m}}_{Q_{k_m}^\nu(X_{D_{m-1}},\cdot)}(X_{k_m},V_{k_m})$ is conditionally $\lambda$-distributed given $U_{D_{m-1}}=u$, so for $\lambda^{D_{m-1}}$-a.e.\ $u$,
    \[
        \left(
        \mathrm{rPIT}^{\preceq_{k_m}}_{Q_{k_m}^\nu(S_{m-1}(u),\cdot)}
        \right)_{\#}
        \left(R_m(S_{m-1}(u),\cdot) \otimes \lambda\right)
        =\lambda.
    \]
    By \cref{thm:randomized_PIT_uniform}, $R_m(S_{m-1}(u),\cdot) = Q_{k_m}^\nu(S_{m-1}(u),\cdot)$ for $\lambda^{D_{m-1}}$-a.e.\ $u$.
    Pushing forward via $S_{m-1}$, we see that the two kernels agree for $\mu_{D_{m-1}}$-a.e.\ history, and thus $\mu_{D_m}=\nu_{D_m}$.

    Finally, if $\bigl(T_{\cI,\nu}^{\vec{\preceq}, \trianglelefteq}\bigr)_{\#}(\mu \otimes \lambda^{\cI})=\lambda^{\cI}$, \cref{prop:grt_projective_consistency} implies that $\bigl(T_{D,\nu_D}^{\vec{\preceq}, \trianglelefteq}\bigr)_{\#}(\mu_D \otimes \lambda^D)=\lambda^D$ for every finite downward-closed $D \subseteq \cI$, whence $\mu_D=\nu_D$ by the finite-prefix result.
    Since any finite $S \subseteq \cI$ is contained in the finite downward-closed set $D(S):=\bigcup_{k \in S}\set{j:j \trianglelefteq k}$, we have $\mu_S=\nu_S$ for all finite $S \subseteq \cI$, and therefore $\mu=\nu$.
\end{proof}

\paragraph{Proof of~\cref{cor:finite_witness_misspecification}.}

\begin{proof}
    Let $\rho:=\bigl(T_{\cI,\nu}^{\vec{\preceq}, \trianglelefteq}\bigr)_{\#}(\mu\otimes\lambda^{\cI})$.
    Since $\mu \neq \nu$, \cref{thm:grt_converse_uniformization} implies $\rho\neq\lambda^{\cI}$, so $(\pi_S)_{\#}\rho\neq\lambda^S$ for some finite $S\subseteq\cI$.
    Taking $D(S):=\bigcup_{k\in S}\set{j:j\trianglelefteq k}$ gives a finite downward-closed set on which $(\pi_{D(S)})_{\#}\rho\neq\lambda^{D(S)}$.
\end{proof}

\paragraph{Proof of~\cref{thm:grt_uniqueness}.}

\begin{proof}
    By construction and \cref{thm:grt_uniformity}, the EGRT $G:=T_{\cI,\nu}^{\vec{\preceq}, \trianglelefteq}$ satisfies properties (i)--(iii).
    Conversely, let $\widetilde T:\cX_{\cI} \times [0,1]^{\cI} \to [0,1]^{\cI}$ satisfy (i)--(iii).
    Enumerating $\cI$ in increasing $\trianglelefteq$-order as $k_1,k_2,\ldots$ (or $k_1,\ldots,k_n$) and setting $D_m:=\{k_1,\ldots,k_m\}$, we prove by induction on $m$ that $\widetilde T_{D_m}=G_{D_m}$ $\nu \otimes \lambda^{\cI}$-a.e.

    For $m=1$, triangularity and local randomization imply that $\widetilde T_{k_1}$ is a function of $(X_{k_1},V_{k_1})$ alone; it is lexicographically nondecreasing and pushes $\nu_{\{k_1\}}\otimes\lambda$ forward to $\lambda$.
    By \cref{lem:one_dimensional_monotone_uniqueness}, $\widetilde T_{k_1}(X,V) = \mathrm{rPIT}_{\nu_{\{k_1\}}}^{\preceq_{k_1}}(X_{k_1},V_{k_1}) = G_{k_1}(X,V)$ a.s.

    Now assume $\widetilde T_{D_{m-1}}=G_{D_{m-1}}$ a.s., and write $k:=k_m$.
    By \cref{lem:conditional_rpit_reconstruction}, there is a measurable reconstruction map $S_{m-1}:[0,1]^{D_{m-1}} \to \cX_{D_{m-1}}$ with $(S_{m-1})_{\#}\lambda^{D_{m-1}}=\nu_{D_{m-1}}$ and $S_{m-1}(G_{D_{m-1}}(X,V))=X_{D_{m-1}}$ a.s.
    Because $\widetilde T$ has product-uniform output and $\widetilde T_{D_{m-1}}=G_{D_{m-1}}$ a.s., the conditional law of $\widetilde T_k$ given $\widetilde T_{D_{m-1}}=u$ is $\lambda$ for $\lambda^{D_{m-1}}$-a.e.\ $u$, and $X_{D_{m-1}}=S_{m-1}(u)$ a.s.\ under this conditioning.
    Moreover, since $X_k \indep V_{D_{m-1}} \mid X_{D_{m-1}}$ and $V_k \indep (X_{\cI},V_{D_{m-1}})$, the conditional law of $(X_k,V_k)$ given $\widetilde T_{D_{m-1}}=u$ is $Q_k^\nu(S_{m-1}(u),dx_k)\lambda(dv_k)$, while local randomization ensures $\widetilde T_k$ does not depend on $V_{D_{m-1}}$.
    For $\lambda^{D_{m-1}}$-a.e.\ $u$, the map $(x_k,v_k) \mapsto \widetilde T_k(S_{m-1}(u),x_k,v_k)$ is lexicographically nondecreasing and pushes $Q_k^\nu(S_{m-1}(u),\cdot) \otimes \lambda$ forward to $\lambda$, so \cref{lem:one_dimensional_monotone_uniqueness} gives
    \[
        \widetilde T_k(S_{m-1}(u),x_k,v_k)
        =
        \mathrm{rPIT}_{Q_k^\nu(S_{m-1}(u),\cdot)}^{\preceq_k}(x_k,v_k)
    \]
    for $Q_k^\nu(S_{m-1}(u),dx_k)\lambda(dv_k)$-a.e.\ $(x_k,v_k)$.
    Integrating over $u \sim \lambda^{D_{m-1}}$ yields $\widetilde T_k(X,V)=G_k(X,V)$ a.s., completing the induction. Taking the union of the null exceptional sets over $m$ gives $\widetilde T=G$ $\nu \otimes \lambda^{\cI}$-a.e.
\end{proof}

\setcounter{section}{2}
\section{Proofs for Section~3}
\label{sec:supp_order_invariance}

\paragraph{Proof of~\cref{lem:map_distributional_equivalence}.}
\begin{proof}
If $T=T'$ $\nu\otimes\lambda^{\cI}$-a.e., then the equality also holds $\mu\otimes\lambda^{\cI}$-a.e.\ for every $\mu\ll\nu$, because $\mu\otimes\lambda^{\cI}\ll\nu\otimes\lambda^{\cI}$; in particular, the two output laws agree under every such $\mu$. This proves the final statement and one direction of the equivalence.

Conversely, suppose that $T$ and $T'$ differ on a set of positive $\nu\otimes\lambda^{\cI}$-measure. Since $\cI$ is countable, there is a coordinate $k$ at which $T_k:=[T]_k$ and $T'_k:=[T']_k$ differ on a set of positive measure. Put $S:=P_{\triangleleft}(k)$, $S':=P_{\triangleleft'}(k)$, and $H:=S\cup S'\cup\{k\}$, and let $Q_k$ and $Q'_k$ be the versions of $Q_{k\mid S}^\nu$ and $Q_{k\mid S'}^\nu$ used to construct $T$ and $T'$ at coordinate $k$. For $z\in\cX_H$ and $u\in[0,1]$, define the fiber maps
\[
    \phi_z(u):=\mathrm{rPIT}^{\preceq_k}_{Q_k(z_S,\cdot)}(z_k,u),
    \qquad
    \phi'_z(u):=\mathrm{rPIT}^{\preceq'_k}_{Q'_k(z_{S'},\cdot)}(z_k,u),
\]
so that $T_k(x,v)=\phi_{x_H}(v_k)$ and $T'_k(x,v)=\phi'_{x_H}(v_k)$. Both maps are jointly measurable in $(z,u)$ by \cref{lem:kernel_valued_rpit_measurability}. Each is nondecreasing and affine in $u$: $\phi_z(u)=a(z)+m(z)u$ with $a(z):=F^{\preceq_k}_{Q_k(z_S,\cdot)}(z_k^-)$ and $m(z):=Q_k(z_S,\{z_k\})\geq0$, and likewise $\phi'_z(u)=a'(z)+m'(z)u$.
Let $K(z,\cdot):=(\phi_z)_{\#}\lambda$ and $K'(z,\cdot):=(\phi'_z)_{\#}\lambda$, which are probability kernels from $\cX_H$ to $[0,1]$. The law $K(z,\cdot)$ is uniform on $[a(z),a(z)+m(z)]$ if $m(z)>0$ and equals $\delta_{a(z)}$ if $m(z)=0$; either way it determines $(a(z),m(z))$, for instance through its mean $a(z)+m(z)/2$ and variance $m(z)^2/12$. The same holds for $K'$.

Under $\nu\otimes\lambda^{\cI}$, the pair $(X_H,V_k)$ has law $\nu_H\otimes\lambda$, so by Fubini's theorem the set
\[
    D_0:=\set{z\in\cX_H:\lambda\set{u:\phi_z(u)\neq\phi'_z(u)}>0}
\]
has $\nu_H(D_0)>0$. For $z\in D_0$, the affine maps $\phi_z$ and $\phi'_z$ are distinct, so $(a(z),m(z))\neq(a'(z),m'(z))$ and hence $K(z,\cdot)\neq K'(z,\cdot)$.
Let $\cA_0$ be a countable algebra generating $\cB([0,1])$. Probability measures that agree on $\cA_0$ agree on $\cB([0,1])$, so $D_0\subseteq\bigcup_{B\in\cA_0}\set{z:K(z,B)\neq K'(z,B)}$, and some $B\in\cA_0$ has $\nu_H\set{z:K(z,B)\neq K'(z,B)}>0$. Replacing $B$ by its complement if necessary, which lies in $\cA_0$ and reverses the sign of $K(z,B)-K'(z,B)$, we may assume that
\[
    D_B:=\set{z\in\cX_H:K(z,B)>K'(z,B)}
\]
has $\nu_H(D_B)>0$.

Let $E:=\set{x\in\cX_{\cI}:x_H\in D_B}$, so that $\nu(E)=\nu_H(D_B)>0$, and let $\mu:=\nu(\cdot\mid E)\ll\nu$. Under $\mu\otimes\lambda^{\cI}$, $V_k$ is uniform and independent of $X_H$, and $X_H$ has law $\nu_H(\cdot\mid D_B)$. Hence
\begin{align*}
    (\mu\otimes\lambda^{\cI})\set{T_k\in B}
    &=\frac{1}{\nu_H(D_B)}\int_{D_B}K(z,B)\,\nu_H(\mathrm dz)\\
    &>\frac{1}{\nu_H(D_B)}\int_{D_B}K'(z,B)\,\nu_H(\mathrm dz)
    =(\mu\otimes\lambda^{\cI})\set{T'_k\in B},
\end{align*}
where the inequality is strict because $K(\cdot,B)-K'(\cdot,B)$ is positive on $D_B$. The $k$-th coordinate marginals of $T_{\#}(\mu\otimes\lambda^{\cI})$ and $T'_{\#}(\mu\otimes\lambda^{\cI})$ therefore differ, and so do the joint output laws.
In particular, whenever the two maps differ, a one-dimensional output marginal already differs under a conditional law $\nu(\cdot\mid E)$ with $\nu(E)>0$.
\end{proof}

\paragraph{Proof of~\cref{prop:countable_extension}.}
\begin{proof}
Let $\preceq$ be a linear extension of $\preceq_p$, and enumerate $\cX=\{x_1,x_2,\ldots\}$ (with the obvious modification if $\cX$ is finite).
We define an order embedding $\phi:(\cX,\preceq)\to(\Q,\leq)$ inductively. Having defined $\phi$ on $\{x_1,\ldots,x_n\}$, set
\[
L_n:=\{\phi(x_i):i \leq n,\ x_i \prec x_{n+1}\},
\qquad
U_n:=\{\phi(x_i):i \leq n,\ x_{n+1} \prec x_i\}.
\]
By the induction hypothesis, $\max L_n<\min U_n$ whenever both sets are nonempty; we may therefore pick $\phi(x_{n+1}) \in \Q$ strictly between $\max L_n$ and $\min U_n$ (treating an empty set as imposing no bound).
By construction, $x \preceq y \Longleftrightarrow \phi(x) \leq \phi(y)$, and since every map on a countable standard Borel space is Borel, $\phi$ is a Borel order embedding.
\end{proof}

\paragraph{Proof of~\cref{thm:admissibility_characterization}.}
The geometric-series construction below adapts \citet[Lemmas~1--3]{hack2022injective} to the Borel setting.
\begin{proof}
(iii)$\Rightarrow$(ii): Let $(f_n)_{n\geq1}$ be a countable point-separating family of monotone Borel functions. Enumerate $\N\times\Q$ as $((n_j,r_j))_{j\geq1}$ and define
\[
h(x):=\sum_{j=1}^\infty3^{-j}\mathbf 1\{f_{n_j}(x)>r_j\}.
\]
Because the series converges uniformly, $h$ is Borel. If $x\neq y$, point separation gives $f_n(x)\neq f_n(y)$ for some $n$, and choosing a rational $r$ between $f_n(x)$ and $f_n(y)$ ensures the two series differ in at least one term. At the first index $j_0$ where they differ, the leading gap $3^{-j_0}$ exceeds the maximum possible tail $\sum_{j>j_0}3^{-j}=3^{-j_0}/2$, so $h(x)\neq h(y)$.
Monotonicity of each $f_n$ implies that if $x \preceq_p y$, then $\mathbf 1\{f_{n_j}(x)>r_j\}\leq\mathbf 1\{f_{n_j}(y)>r_j\}$ for all $j$, giving $h(x)\leq h(y)$; if $x\prec_p y$, then $h(x)<h(y)$ by injectivity.
(Conversely, if $h:\cX\to\R$ satisfies (ii), the singleton family $\{h\}$ immediately gives (iii).)

(ii)$\Rightarrow$(i): Given $h$ as in (ii), define $x\preceq_h y \Longleftrightarrow h(x)\leq h(y)$. Since $h$ is injective, $\preceq_h$ is a total order extending $\preceq_p$ with Borel order embedding $h$, which also establishes the final claim of the theorem.

(i)$\Rightarrow$(ii): If $\preceq\in\cE(\preceq_p)$ has Borel order embedding $h_\preceq$, then $x\prec_p y$ implies $x\prec y$ and thus $h_\preceq(x)<h_\preceq(y)$.
\end{proof}

\begin{lemma}[Measurability of conditional ambiguity]
\label{lem:conditional_ambiguity_measurable}
In the setting of~\cref{defn:conditional_ambiguity}, the map
\[
\Delta_k(x_{\trianglelefteq k})
=Q_k^\nu(x_{\triangleleft k},I_{p,k}(x_k))
\]
is $\cB_{\trianglelefteq k}$-measurable.
\end{lemma}

\begin{proof}
Since the graph $G_k:=\{(x,y):x \preceq_{p,k}y\}$ is Borel by assumption, the incomparability relation
$J_k:=(\cX_k \times \cX_k) \setminus (G_k \cup G_k^\top)$ is also Borel, with $x$-section $I_{p,k}(x)$.
Viewing $Q_k^\nu(x_{\triangleleft k},\mathrm dy)$ as a kernel $\widetilde Q_k^\nu((x_{\triangleleft k},x_k),\mathrm dy)$ on $\cX_{\trianglelefteq k}\times\cB_k$, standard measurability of kernel integrals shows that
\[
(x_{\triangleleft k},x_k)
\longmapsto
\int\mathbf 1_{J_k}(x_k,y)\,
Q_k^\nu(x_{\triangleleft k},\mathrm dy)
=\Delta_k(x_{\trianglelefteq k})
\]
is $\cB_{\trianglelefteq k}$-measurable.
\end{proof}

\begin{lemma}[CDF separation on countable spaces]
\label{lem:countable_cdf_separating}
Let $(\cX,\cB)$ be a countable standard Borel space with partial order
$\preceq_p$. Then $\cE(\preceq_p)$ is CDF-separating. A countable
separating family can be chosen to contain, for every ordered incomparable
pair $(a,b)$, an extension with $a \prec b$. If $\cX$ is finite, the
finite family of all linear extensions may be used.
\end{lemma}

\begin{proof}
For each ordered incomparable pair $(a,b)$, adjoining $a \prec b$ and taking the transitive closure introduces no cycles, so Szpilrajn's theorem yields a linear extension $\preceq_{a,b}$, which is Borel-embeddable by \cref{prop:countable_extension}. Because $\cX$ is countable, the collection of such extensions $\preceq_{a,b}$ (augmented by a single extension if $\preceq_p$ is already total) is countable.

If $\eta\{x:\eta(I_p(x))>0\}>0$, countability of $\cX$ ensures there exist incomparable $x,y \in \cX$ with $\eta(\{x\})>0$ and $\eta(\{y\})>0$. Setting $\preceq^+ := \preceq_{x,y}$ and $\preceq^- := \preceq_{y,x}$, we have
\[
    F_\eta^{\preceq^+}(y)-F_\eta^{\preceq^+}(x)
    \geq \eta(\{y\})>0,
    \qquad
    F_\eta^{\preceq^-}(y)-F_\eta^{\preceq^-}(x)
    \leq -\eta(\{x\})<0,
\]
so $F_\eta^{\preceq^+}$ and $F_\eta^{\preceq^-}$ cannot agree at both $x$ and $y$. Because the point masses $\eta(\{x\})$ and $\eta(\{y\})$ do not depend on the order, the strict-left CDFs likewise differ on a set of positive $\eta$-measure. When $\cX$ is finite, the set of all linear extensions is finite and contains $\preceq^+$ and $\preceq^-$.
\end{proof}

\paragraph{Proof of~\cref{thm:countable_representation_cdf_separating}.}
While \citet[Proposition~6]{hack2022injective} relate countable multi-utilities to countable families of injective monotones, we must show here that Borel extensions separate CDFs in the presence of positive incomparability mass (even for atomless laws) and that this probabilistic separation recovers the order representation.
\begin{proof}
(i)$\Rightarrow$(ii): By \cref{eq:countable_order_representation}, the graph of $\preceq_p$ is $\bigcap_{n\geq 1}\{(x,y):f_n(x)\leq f_n(y)\}$ for a sequence of $\preceq_p$-monotone Borel functions $(f_n)_{n\geq 1}$. By antisymmetry of $\preceq_p$, $(f_n)_{n\geq 1}$ separates points, so \cref{thm:admissibility_characterization} (composed with a strictly increasing bijection $\R \to (0,1)$) yields a Borel injection $h_0:\cX\to(0,1)$ satisfying $x\prec_p y \Rightarrow h_0(x)<h_0(y)$.

For each $n\geq1$ and $q\in\Q$, define $h_{n,q}:=\mathbf1\{f_n>q\}+h_0$. Since $h_0(\cX)\subseteq(0,1)$, $h_{n,q}$ maps $\{f_n\leq q\}$ into $(0,1)$ and $\{f_n>q\}$ into $(1,2)$, making $h_{n,q}$ Borel and injective; furthermore, $x\prec_p y$ implies $h_{n,q}(x)<h_{n,q}(y)$. Thus its pullback order $\preceq_{n,q}$ belongs to $\cE(\preceq_p)$, and $\cE_0:=\{\preceq_{n,q}:n\geq1,\ q\in\Q\}$ is a countable subfamily of $\cE(\preceq_p)$.

If $x,y\in\cX$ are incomparable, \cref{eq:countable_order_representation} gives $n,m$ with $f_n(x)>f_n(y)$ and $f_m(y)>f_m(x)$. Picking rationals $q, q'$ strictly between these pairs of values, we see that $\preceq_{n,q}$ and $\preceq_{m,q'}$ order $x$ and $y$ in opposite directions. Hence the incomparability set $J:=\{(x,y):y\in I_p(x)\}$ is covered by the countable union
\[
    J\subseteq
    \bigcup_{\preceq,\preceq'\in\cE_0}
    R_{\preceq,\preceq'},
    \qquad
    R_{\preceq,\preceq'}
    :=\{(x,y):x\prec y,\ y\prec' x\}.
\]
Now let $\eta$ be any probability measure with $\eta\{x:\eta(I_p(x))>0\}>0$. By Fubini's theorem, $(\eta\otimes\eta)(J)>0$, so $(\eta\otimes\eta)(R_{\preceq,\preceq'})>0$ for some pair $\preceq,\preceq'\in\cE_0$.
Suppose for contradiction that $F_\eta^\preceq(\cdot^-) = F_\eta^{\preceq'}(\cdot^-)$ $\eta$-a.e.; since singleton masses are order-independent, we then have $\mathrm{rPIT}_\eta^\preceq(x,v)=\mathrm{rPIT}_\eta^{\preceq'}(x,v)$ for all $v\in[0,1]$ and $\eta$-a.e.\ $x$.
Let $(X,Y,V,W) \sim \eta\otimes\eta\otimes\lambda\otimes\lambda$, and put
\[
    U:=\mathrm{rPIT}_\eta^\preceq(X,V),
    \qquad
    Z:=\mathrm{rPIT}_\eta^\preceq(Y,W).
\]
By \cref{thm:randomized_PIT_uniform}, $U$ and $Z$ are independent $\mathrm{Unif}(0,1)$ random variables. Yet whenever $x\prec y$,
\[
    \mathrm{rPIT}_\eta^\preceq(x,v)
    \leq F_\eta^\preceq(x)
    \leq F_\eta^\preceq(y^-)
    \leq \mathrm{rPIT}_\eta^\preceq(y,w)
    \qquad\text{for all }v,w\in[0,1],
\]
which implies $U\leq Z$ on $\{(X,Y)\in R_{\preceq,\preceq'}\}$. Applying the same inequality to $\preceq'$ (where $Y\prec' X$) and using the a.s.\ agreement of the two rPITs gives $Z\leq U$ a.s.\ on $\{(X,Y)\in R_{\preceq,\preceq'}\}$. Thus $U=Z$ on an event of positive probability, which is impossible for independent continuous random variables. Hence $\cE_0$ witnesses CDF separation for every probability measure $\eta$.

(ii)$\Rightarrow$(i): Choose a countable family $\cE_0 \subseteq \cE(\preceq_p)$ witnessing CDF separation for all probability laws. Since $\cE(\preceq_p)\neq\varnothing$, we may adjoin an admissible extension if necessary to make $\cE_0$ nonempty. For each $\preceq\in\cE_0$, choose a Borel order embedding $f_{\preceq}:\cX\to\R$.
For any incomparable pair $x,y$, apply CDF separation to $\eta=(\delta_x+\delta_y)/2$. Since $\eta(I_p(x))=\eta(I_p(y))=1/2>0$, the strict-left CDFs cannot take identical values on $\{x,y\}$ across all $\preceq\in\cE_0$, meaning $\cE_0$ must contain extensions ordering $(x,y)$ in both directions; in particular, $f_{\preceq}(x)>f_{\preceq}(y)$ for some $\preceq\in\cE_0$ and $f_{\preceq'}(x)<f_{\preceq'}(y)$ for another.
Consequently,
\[
    x\preceq_p y
    \quad\Longleftrightarrow\quad
    f_{\preceq}(x)\leq f_{\preceq}(y)
    \ \text{for every }\preceq\in\cE_0,
\]
and enumerating $\{f_{\preceq}:\preceq\in\cE_0\}$ gives (i).
\end{proof}

\paragraph{Proof of~\cref{prop:zero_ambiguity_sufficient}.}
\begin{proof}
Fix $k\in\cI$, two extensions $\preceq,\preceq'$ of $\preceq_{p,k}$, a probability measure $\eta$ on $\cX_k$, and $x\in\cX_k$. Setting $L:=\{y:y \prec x\}$ and $L':=\{y:y \prec' x\}$, note that any $y\in L\mathbin\triangle L'$ is ordered oppositely relative to $x$ by $\preceq$ and $\preceq'$, so $L\mathbin\triangle L'\subseteq I_{p,k}(x)$. Because $\eta(\{x\})$ does not depend on the order, we have for all $v\in[0,1]$,
\begin{equation}
\left|
\mathrm{rPIT}_{\eta}^{\preceq}(x,v)
-
\mathrm{rPIT}_{\eta}^{\preceq'}(x,v)
\right|
=
\left|\eta(L)-\eta(L')\right|
\leq
\eta(L\mathbin\triangle L')
\leq
\eta\bigl(I_{p,k}(x)\bigr).
\label{eq:supp_component_delta_pointwise}
\end{equation}
Taking $\eta=Q_k^\nu(x_{\triangleleft k},\cdot)$ and $x=x_k$ in \cref{eq:supp_component_delta_pointwise} yields the pointwise bound in the proposition.
When $\Delta_k=0$ $\nu_{\trianglelefteq k}$-a.e.\ for all $k\in\cI$, the two transforms agree coordinatewise $\nu\otimes\lambda^{\cI}$-a.e., and hence globally $\nu\otimes\lambda^{\cI}$-a.e.\ since $\cI$ is countable; by absolute continuity, this a.e.\ equality and equality of the pushforward laws carry over to every $\mu\ll\nu$.
\end{proof}

\paragraph{Proof of~\cref{thm:component_order_exact_invariance}.}
\begin{proof}
\Cref{prop:zero_ambiguity_sufficient} gives (i)$\Rightarrow$(ii), and \cref{lem:map_distributional_equivalence} gives (ii)$\Leftrightarrow$(iii) together with the final statement of the theorem. Neither uses CDF separation.

Now assume that $\cE(\preceq_{p,k})$ is CDF-separating for $\cH_k^\nu$ for every $k\in\cI$. To prove (ii)$\Rightarrow$(i), suppose by contraposition that (i) fails at coordinate $k$. Writing $h=x_{\triangleleft k}$ and $\eta_h=Q_k^\nu(h,\cdot)$, we have
\[
\int
\eta_h\{x:\eta_h(I_{p,k}(x))>0\}\,
\nu_{\triangleleft k}(\mathrm dh)>0,
\]
so $\{h:\eta_h\{x:\eta_h(I_{p,k}(x))>0\}>0\}$ has positive $\nu_{\triangleleft k}$-measure.
Let $\cE_0(\preceq_{p,k})$ be a countable subfamily witnessing CDF separation for $\cH_k^\nu$ in~\cref{eq:reference_kernel_family}. Because $\cE_0(\preceq_{p,k})$ is countable, there is a fixed pair $\preceq_k, \preceq_k' \in \cE_0(\preceq_{p,k})$ whose strict-left CDFs under $\eta_h$ differ on a set of positive $\eta_h$-measure for a set of histories $h$ of positive $\nu_{\triangleleft k}$-measure.
By \cref{lem:kernel_valued_rpit_measurability} and Fubini's theorem, the set
\[
D:=\{(h,x_k):
F_{\eta_h}^{\preceq_k}(x_k^-)
\neq
F_{\eta_h}^{\preceq_k'}(x_k^-)
\}
\]
is measurable with $\nu_{\trianglelefteq k}(D)>0$. Extend $\preceq_k$ and $\preceq_k'$ to full component-order sequences $\vec\preceq$ and $\vec\preceq'$ that coincide at all coordinates other than $k$.

Both transforms use the same kernel $Q_k^\nu$ at coordinate $k$, and the atom $\eta_h(\{x_k\})$ does not depend on the order. Hence, for every $(x,v)\in\cX_{\cI}\times[0,1]^{\cI}$ with $x_{\trianglelefteq k}=(h,x_k)$,
\[
[T_{\cI,\nu}^{\vec\preceq}(x,v)]_k-[T_{\cI,\nu}^{\vec\preceq'}(x,v)]_k
=F_{\eta_h}^{\preceq_k}(x_k^-)-F_{\eta_h}^{\preceq_k'}(x_k^-),
\]
which is nonzero whenever $(h,x_k)\in D$, whatever the value of $v$. The two transforms therefore differ on $\set{(x,v):x_{\trianglelefteq k}\in D}$, which has $\nu\otimes\lambda^{\cI}$-measure $\nu_{\trianglelefteq k}(D)>0$, contradicting (ii).
\end{proof}

\begin{lemma}[Existence of finite-predecessor linear extensions]
\label{lem:index_linear_extensions_exist}
Let $\trianglelefteq_p$ be a partial order on a countable set $\cI$ such that every $\mathrm{anc}_p(k)$ is finite. Then $\mathcal L(\trianglelefteq_p)$ is nonempty.
\end{lemma}

\begin{proof}
Fix an enumeration $\cI=\{i_1,i_2,\ldots\}$ and recursively construct a sequence by selecting at each step the $\trianglelefteq_p$-minimal remaining element of $\cI$ with the smallest enumeration index. (A minimal remaining element always exists, as an infinite descending chain would yield infinitely many ancestors.)
By induction on $|\mathrm{anc}_p(k)|$, every $k\in\cI$ is eventually selected: once all elements of $\mathrm{anc}_p(k)$ have been chosen, $k$ becomes minimal among the remaining elements and can be preceded by at most finitely many indices with smaller enumeration index. Ordering $\cI$ by selection time gives an extension in $\mathcal L(\trianglelefteq_p)$.
\end{proof}

\begin{lemma}[Control of predecessor sets in index extensions]
\label{lem:flexible_index_linear_extensions}
Under the assumptions of~\cref{lem:index_linear_extensions_exist}, fix $k \in \cI$.
\begin{enumerate}[label=(\alph*)]
    \item There exists $\trianglelefteq_k^0 \in \mathcal L(\trianglelefteq_p)$ such that
    \[
        P_{\triangleleft_k^0}(k)=\mathrm{anc}_p(k).
    \]
    \item For every finite $B \subseteq \mathrm{inc}_p(k)$, there exists
    $\trianglelefteq^{k,B} \in \mathcal L(\trianglelefteq_p)$ such that
    \[
        B \subseteq P_{\triangleleft^{k,B}}(k).
    \]
\end{enumerate}
\end{lemma}

\begin{proof}
(a) Adjoin to $\trianglelefteq_p$ the relations $k\triangleleft j$ for all $j\in\mathrm{inc}_p(k)$ and take the transitive closure. Because every added relation exits $k$ and no $j\in\mathrm{inc}_p(k)$ is an ancestor of $k$, no cycle can arise. Under this enlarged partial order, the ancestor set of any index is contained in its original ancestor set united with $\mathrm{anc}_p(k)\cup\{k\}$, which is finite. Applying \cref{lem:index_linear_extensions_exist} to the enlarged partial order yields $\trianglelefteq_k^0 \in \mathcal L(\trianglelefteq_p)$ in which every $j \notin \mathrm{anc}_p(k)\cup\{k\}$ follows $k$, so $P_{\triangleleft_k^0}(k)=\mathrm{anc}_p(k)$.

(b) Similarly, adjoining $b\triangleleft k$ for each $b\in B$ and taking the transitive closure creates no cycles, and increases any ancestor set by at most the finite set $\bigcup_{b \in B}(\mathrm{anc}_p(b) \cup \{b\})$. Any extension in $\mathcal L$ of this enlarged partial order (via \cref{lem:index_linear_extensions_exist}) belongs to $\mathcal L(\trianglelefteq_p)$ and satisfies $B \subseteq P_{\triangleleft^{k,B}}(k)$.
\end{proof}

% Keep the proof heading with its opening paragraph.
\begingroup
\makeatletter\@beginparpenalty=10000\makeatother
\paragraph{Proof of~\cref{thm:exact_index_order_invariance}.}
\begin{proof}
Fix $k \in \cI$ and put $A:=\mathrm{anc}_p(k)$ and $C:=\mathrm{inc}_p(k)$. For any $\trianglelefteq \in \mathcal L(\trianglelefteq_p)$, $S:=P_{\triangleleft}(k)$ satisfies $A \subseteq S \subseteq A \cup C$.

(i)$\Leftrightarrow$(ii): If (i) holds, decomposition of conditional independence yields $X_k \indep_\nu X_{S \setminus A} \mid X_A$, which is equivalent to the kernel identity in (ii).
Conversely, assume (ii) and let $B\subseteq C$ be finite. By \cref{lem:flexible_index_linear_extensions}(b), there exists $\trianglelefteq^{k,B} \in \mathcal L(\trianglelefteq_p)$ with $B\subseteq P_{\triangleleft^{k,B}}(k) \setminus A$, so (ii) gives $X_k \indep_\nu X_{P_{\triangleleft^{k,B}}(k) \setminus A} \mid X_A$ and hence $X_k\indep_\nu X_B\mid X_A$. Since $B\subseteq C$ was an arbitrary finite subset and $C$ is countable, a monotone-class argument yields $X_k\indep_\nu X_C\mid X_A$.

(ii)$\Leftrightarrow$(iii): Under (ii), for every $\trianglelefteq \in \mathcal L(\trianglelefteq_p)$ and $k\in\cI$,
\[
    [T_{\cI,\nu}^{\trianglelefteq}(X,V)]_k
    =
    \mathrm{rPIT}^{\preceq_k}_{
        Q_{k \mid A}^\nu(X_A,\cdot)}(X_k,V_k)
    \qquad \nu \otimes \lambda^{\cI}\text{-a.e.}
\]
Since the right-hand side is independent of $\trianglelefteq$ and $\cI$ is countable, (iii) follows.
Conversely, suppose (iii) holds and fix $\trianglelefteq \in \mathcal L(\trianglelefteq_p)$ with $S=P_{\triangleleft}(k)$. Comparing $\trianglelefteq$ against the extension $\trianglelefteq_k^0$ from \cref{lem:flexible_index_linear_extensions}(a) (for which $P_{\triangleleft_k^0}(k)=A$) gives
\begin{align*}
&\mathrm{rPIT}^{\preceq_k}_{
    Q_{k \mid S}^\nu(X_S,\cdot)}(X_k,V_k)
\\
&\qquad=
\mathrm{rPIT}^{\preceq_k}_{
    Q_{k \mid A}^\nu(X_A,\cdot)}(X_k,V_k)
\qquad \nu \otimes \lambda^{\cI}\text{-a.e.}
\end{align*}
Disintegrating over $X_S=s$, we see that for $\nu_S$-a.e.\ $s$, the two conditional rPITs agree $Q_{k\mid S}^\nu(s,\mathrm dx)\lambda(\mathrm dv)$-a.e. Because the left-hand rPIT pushes $Q_{k\mid S}^\nu(s,\cdot)\otimes\lambda$ forward to $\lambda$, so does the right-hand rPIT, and \cref{thm:randomized_PIT_uniform} therefore implies $Q_{k \mid S}^\nu(s,\cdot) = Q_{k \mid A}^\nu(s_A,\cdot)$ for $\nu_S$-a.e.\ $s$, which is (ii).

(iii)$\Leftrightarrow$(iv): This is \cref{lem:map_distributional_equivalence} applied to $T_{\cI,\nu}^{\trianglelefteq}$ and $T_{\cI,\nu}^{\trianglelefteq'}$ for each pair $\trianglelefteq,\trianglelefteq'\in\mathcal L(\trianglelefteq_p)$. The same lemma gives the final statement of the theorem.
\end{proof}
\endgroup

\paragraph{Proof of~\cref{cor:dag_ancestral_screening}.}
\begin{proof}
Fix $k\in\cI$ and put $A:=\mathrm{anc}_p(k)$, $C:=\mathrm{inc}_p(k)$, and $P:=\mathrm{Pa}_G(k)$, so that $P\subseteq A$ and $\mathrm{ND}_G(k)=A\cup C$ (disjoint). By the local Markov property,
\[
    Q_{k\mid A\cup C}^\nu(X_{A\cup C},\cdot)
    =
    Q_{k\mid P}^\nu(X_P,\cdot)
    \qquad\nu\text{-a.s.}
\]
Taking conditional expectations given $X_A$ and using $P\subseteq A$, we also have $Q_{k\mid A}^\nu(X_A,\cdot) = Q_{k\mid P}^\nu(X_P,\cdot)$ $\nu$-a.s., whence $X_k\indep_\nu X_C\mid X_A$. Since every $k\in\cI$ has a finite ancestor set by assumption, \cref{thm:exact_index_order_invariance} applies.
\end{proof}

\setcounter{section}{3}
\section{Proofs for Section~4}
\label{sec:supp_statistical_equivalence}

\paragraph{Proof of~\cref{lem:egrt_output_reconstruction}.}
\begin{proof}
Following the recursive quantile-inverse construction of \citet[Theorem~3]{ruschendorf1993regression} \citep[see also][Section~3]{ruschendorf2009distributional}, enumerate $\cI$ in increasing $o$-index order as $k_1 \triangleleft k_2 \triangleleft \cdots$ (or $k_1\triangleleft\cdots\triangleleft k_n$) and write $D_m:=\{k_1,\ldots,k_m\}$.
Applying \cref{lem:conditional_rpit_reconstruction} to the marginal at $k_1$ gives a measurable map $g_1:[0,1] \to \cX_{k_1}$ reconstructing $X_{k_1}$, and applying \cref{lem:conditional_rpit_reconstruction} inductively to the kernel $Q_{k_m}^\nu$ on $\cX_{D_{m-1}} \times \cB_{k_m}$ gives measurable maps $g_m:\cX_{D_{m-1}} \times [0,1] \to \cX_{k_m}$ for $m\geq 2$.
We then define $G_o:[0,1]^{\cI}\to\cX_{\cI}$ coordinatewise by
\[
    [G_o(u)]_{k_1}:=g_1(u_{k_1}),
    \qquad
    [G_o(u)]_{k_m}
    :=
    g_m\bigl(G_o(u)_{D_{m-1}},u_{k_m}\bigr)
    \quad (m\geq 2).
\]
By induction, each coordinate map is measurable, and hence $G_o$ is measurable into the product space $(\cX_{\cI},\cB_{\cI})$; note that $G_o$ depends only on $\nu$ and $o$.

Under $\nu \otimes \lambda^{\cI}$, induction on $m$ via \cref{lem:conditional_rpit_reconstruction} shows that $G_o(T_o(X,V))_{D_m}=X_{D_m}$ a.s.\ for all $m$, and taking a countable union over $m$ gives \cref{eq:main_egrt_reconstruction} under $\nu \otimes \lambda^{\cI}$.
Now let $\mu$ be locally dominated by $\nu$. For each $m$, both
$T_o(x,v)_{D_m}$ and $G_o(u)_{D_m}$ depend only on their respective
inputs on $D_m$, by triangularity and the recursive construction.
Thus failure of reconstruction on $D_m$ is determined by a
measurable subset of $\cX_{D_m}\times[0,1]^{D_m}$ that is null under
$\nu_{D_m}\otimes\lambda^{D_m}$. Since
$\mu_{D_m}\otimes\lambda^{D_m}\ll
\nu_{D_m}\otimes\lambda^{D_m}$, it is also null under
$\mu_{D_m}\otimes\lambda^{D_m}$. Taking a countable union over
$m$ proves \cref{eq:main_egrt_reconstruction} under
$\mu\otimes\lambda^{\cI}$. Pushing forward
$T_o(X,V)\sim\rho_o^\mu$ by $G_o$ then yields
$(G_o)_{\#}\rho_o^\mu=\mu$.
\end{proof}

\paragraph{Proof of~\cref{thm:exact_testing_equivalence}.}
\begin{proof}
For every $\mu\in\mathcal M$, we have $\rho_o^\mu=\mu\mathsf P_o$ by definition and $\rho_o^\mu\mathsf R_o=\mu$ by \cref{eq:reverse_kernel_recovers_input_law}, where neither Markov kernel $\mathsf P_o$ nor $\mathsf R_o$ depends on $\mu$. By the randomization criterion \citep{blackwell1953equivalent,torgersen1991comparison}, the original and transformed experiments are equivalent.
Similarly, for any $o,o' \in \fO$,
\[
    \rho_o^\mu(\mathsf R_o\mathsf P_{o'})
    =
    (\rho_o^\mu\mathsf R_o)\mathsf P_{o'}
    =
    \mu\mathsf P_{o'}
    =
    \rho_{o'}^\mu
\]
and $\rho_{o'}^\mu(\mathsf R_{o'}\mathsf P_o)=\rho_o^\mu$, so the $\mu$-independent composite kernels $\mathsf R_o\mathsf P_{o'}$ and $\mathsf R_{o'}\mathsf P_o$ establish equivalence between any two transformed experiments.
\end{proof}

\paragraph{Proof of~\cref{cor:pairwise_test_transport}.}
\begin{proof}
Because $\mathsf P_o$ is a Markov kernel and $G_{o'}$ is measurable, $\psi'=(\mathsf P_o\psi)\circ G_{o'}$ is a $[0,1]$-valued measurable test, and \cref{lem:egrt_output_reconstruction} ($(G_{o'})_{\#}\rho_{o'}^\mu=\mu$) gives
\[
    \int\psi'\,\mathrm d\rho_{o'}^\mu
    =\int\mathsf P_o\psi\,\mathrm d\mu
    =\int\psi\,\mathrm d(\mu\mathsf P_o)
    =\int\psi\,\mathrm d\rho_o^\mu
    \qquad\text{for all locally dominated }\mu.
\]
Moreover, if $U'\sim\rho_{o'}^\mu$ and $V\sim\lambda^{\cI}$ is drawn independently, then $G_{o'}(U')\sim\mu$ and hence $T_o(G_{o'}(U'),V)\sim\rho_o^\mu$; conditional on $U'=u'$, the randomized draw $T_o(G_{o'}(u'),V)$ rejects under $\psi$ with probability $\psi'(u')$ by \cref{eq:pairwise_test_transport}.
\end{proof}

% Keep the proof heading with its opening paragraph.
\begingroup
\makeatletter\@beginparpenalty=10000\makeatother
\paragraph{Proof of~\cref{prop:egrt_output_likelihood_ratio}.}
\begin{proof}
For any measurable $h:[0,1]^{\cI}\to[0,\infty]$, applying \cref{eq:main_egrt_reconstruction} under $\nu\otimes\lambda^{\cI}$ and using $\rho_o^\nu=\lambda^{\cI}$ gives
\[
\begin{split}
    \int h(u)\,\rho_o^\mu(\mathrm du)
    &=\int h(T_o(x,v))R(x)\,
        (\nu\otimes\lambda^{\cI})(\mathrm d(x,v))\\
    &=\int h(T_o(x,v))R(G_o(T_o(x,v)))\,
        (\nu\otimes\lambda^{\cI})(\mathrm d(x,v))\\
    &=\int h(u)R(G_o(u))\,\lambda^{\cI}(\mathrm du),
\end{split}
\]
from which $\rho_o^\mu\ll\lambda^{\cI}$ and \cref{eq:egrt_output_likelihood_ratio} follow by taking $h$ to be indicator functions.
\end{proof}
\endgroup

\paragraph{Proof of~\cref{cor:egrt_divergence_preservation}.}
\begin{proof}
Because $\mathsf P_o$ maps $(\mu,\nu)$ to $(\rho_o^\mu,\lambda^{\cI})$ and $\mathsf R_o$ maps $(\rho_o^\mu,\lambda^{\cI})$ back to $(\mu,\nu)$ by \cref{eq:reverse_kernel_recovers_input_law}, applying the data-processing inequality \citep[Theorem~3.1]{liese2012divergences} in both directions gives
\[
    D_f(\mu\|\nu)
    \leq D_f\!\left(\rho_o^\mu\middle\|\lambda^{\cI}\right)
    \leq D_f(\mu\|\nu).
\]
\end{proof}

\paragraph{Proof of~\cref{cor:egrt_finite_prefix_divergence}.}
\begin{proof}
Write $o=(\vec{\preceq},\trianglelefteq)$ and let $T_D:=T_{D,\nu_D}^{\vec{\preceq},\trianglelefteq}$ be the finite-prefix EGRT of $\nu_D$ defined prior to the proof of \cref{thm:grt_uniformity}.
By \cref{prop:grt_projective_consistency}, $\pi_D^{[0,1]}\circ T_o=T_D\circ(\pi_D^{\cX}\times\pi_D^{[0,1]})$ pointwise, so
\[
    \rho_{o,D}^\mu
    =\bigl(\pi_D^{[0,1]}\bigr)_{\#}\rho_o^\mu
    =(T_D)_{\#}(\mu_D\otimes\lambda^D).
\]
Since local domination gives $\mu_D\ll\nu_D$, applying \cref{cor:egrt_divergence_preservation} to $T_D$ on the finite index set $D$ yields \cref{eq:egrt_finite_prefix_divergence}.
\end{proof}

\paragraph{Proof of~\cref{cor:egrt_finite_prefix_divergence_limit}.}
\begin{proof}
By \cref{cor:egrt_finite_prefix_divergence}, $D_f(\rho_{o,D_m}^\mu\|\lambda^{D_m}) = D_f(\mu_{D_m}\|\nu_{D_m})$. Since $(\mu_{D_m},\nu_{D_m})$ is a coordinate projection of both $(\mu_{D_{m+1}},\nu_{D_{m+1}})$ and $(\mu,\nu)$, the data-processing inequality \citep[Theorem~3.1]{liese2012divergences} implies that $D_f(\mu_{D_m}\|\nu_{D_m})$ is nondecreasing in $m$ and bounded above by $D_f(\mu\|\nu)$.

For the reverse inequality, work on $(\cX_{\cI},\cB_{\cI},\nu)$, put $\mathcal F_m:=\sigma(\pi_{D_m})$, and let $R:=\mathrm d\mu/\mathrm d\nu$. Then $R_m:=\mathbb E_\nu[R\mid\mathcal F_m]$ is a version of $(\mathrm d\mu_{D_m}/\mathrm d\nu_{D_m})(X_{D_m})$, and since $D_m\uparrow\cI$ implies $\bigvee_m \mathcal F_m = \cB_{\cI}$, L\'evy's upward martingale convergence theorem \citep[Chapter~9]{kallenberg2021foundations} gives $R_m\to R$ $\nu$-a.s.
Choosing a subgradient $a\in\R$ of $f$ at $1$, the nonnegative convex function $\widetilde f(t):=f(t)-a(t-1)$ satisfies
\[
    D_f(\mu_{D_m}\|\nu_{D_m})=\mathbb E_\nu[\widetilde f(R_m)],
    \qquad
    D_f(\mu\|\nu)=\mathbb E_\nu[\widetilde f(R)].
\]
The right-limit extension at zero makes $\widetilde f$ lower semicontinuous on $[0,\infty)$, including when $\widetilde f(0)=+\infty$. Thus $\widetilde f(R)\leq\liminf_m\widetilde f(R_m)$ $\nu$-a.s., and Fatou's lemma yields $\liminf_{m\to\infty}D_f(\mu_{D_m}\|\nu_{D_m}) \geq D_f(\mu\|\nu)$.
\end{proof}

\end{document}